\documentclass[11pt]{amsart}

\usepackage{amssymb,amsmath}
\usepackage{graphics,epsfig,psfrag}
\usepackage{a4wide,mathrsfs}

\numberwithin{equation}{section}

\newcommand{\calA}{\mathcal{A}}

\newcommand{\mC}{\mathbb{C}}
\newcommand{\mD}{\mathbb{D}}

\newcommand{\mN}{\mathbb{N}}

\newcommand{\mR}{\mathbb{R}}

\newcommand{\mT}{\mathbb{T}}

\newcommand{\mZ}{\mathbb{Z}}

\newcommand{\p}{\partial}
\newcommand{\pb}{\overline\partial}
\newcommand{\C}{\mathbb{C}}
\newcommand{\T}{\mathbb{T}}

\newcommand{\e}{\varepsilon}
\newcommand{\f}{\varphi}

\newcommand{\wt}{\widetilde}
\newcommand{\cA}{\mathcal{A}}
\newcommand{\fdot}{\,\cdot\,}

\newcommand{\ci}[1]{_{ {}_{\scriptstyle #1}}}

\newcommand{\nm}{\,\rule[-.6ex]{.13em}{2.3ex}\,}

\newtheorem{theorem}{Theorem}[section]
\newtheorem{lemma}[theorem]{Lemma}

\newtheorem{proposition}[theorem]{Proposition}

\theoremstyle{definition}
\newtheorem{rem}[theorem]{Remark}

\theoremstyle{definition}
\newtheorem{df}[theorem]{Definition}

\newenvironment{entry}
{\begin{list}{X}
  {
      \setlength{\labelwidth}{55pt}
      \setlength{\leftmargin}{\labelwidth}
      \addtolength{\leftmargin}{\labelsep}
      \setlength{\itemsep}{.4pc}
   }
}
{\end{list}}

\begin{document}

\title[Estimates in corona theorems for some subalgebras
  of $H^{\infty}$]{Estimates in corona theorems \\
  for some subalgebras of $H^{\infty}$}

\author{Amol Sasane}
\address{
Department of Mathematics,
London School of Economics,
Houghton Street,
London WC2A 2AE,
United Kingdom.
}
\email{A.J.Sasane@lse.ac.uk}
\author{Sergei Treil}
\thanks{The work of S.~Treil was supported by the National Science Foundation under Grant  DMS-0501065}
\address{
Mathematics Department,
Brown University,
151 Thayer Street/Box 1917,
Providence, RI  02912, U.S.A.
}
\email{treil@math.brown.edu}


%\date{25 May, 2006}


\begin{abstract}
If $n$ is a nonnegative integer, then denote by
$\partial^{-n}H^{\infty}$ the space of all complex valued functions
$f$ defined on $\mathbb{D}$ such that $f,f^{(1)},f^{(2)},\dots,f^{(n)}$
belong to $H^{\infty}$, with the norm
\[
\|f\| =
\sum_{j=0}^{n}\frac{1}{j!}\|f^{(j)}\|_{\infty}.
\]
We prove bounds on the solution in the corona problem for
$\partial^{-n}H^{\infty}$. As corollaries, we obtain estimates
in the corona theorem also for some other subalgebras of the
Hardy space $H^{\infty}$.
\end{abstract}

\maketitle

\tableofcontents

\newpage

\section*{Notation}

\begin{entry}

\item[$:=$]   equal by definition;

\item[$\mC$]  the complex plane;

\item[$\mD$]  the unit disk, $\mD:=\{z \in \mC\;|\; |z|<1\}$;

\item[$\overline{\mD}$] the closed unit disk,
$\overline{\mD}:=\{z \in \mC\;|\;|z|<1\}$;

\item[$\mT$]  the unit circle, $\mT:=\partial \mD
=\{z \in \mC\;|\;|z|=1\}$;

\item[$dm$]  normalized Lebesgue measure on $\mT$,
$m(\mT)=1$;

\item[$\partial, \overline{\partial}$]   derivatives
with respect to $z$ and $\overline{z}$ respectively:
$\partial
\!:=\!
\frac{1}{2}(\frac{\partial}{\partial x}
+i\frac{\partial}{\partial y})$,
$\overline{\partial}
\!:=\!
\frac{1}{2}(\frac{\partial}{\partial x}
-i\frac{\partial}{\partial y})$;

\item[$\Delta$]  Laplacian,
$\Delta:=4\partial \overline{\partial}$;

\item[$\nm \cdot \nm$, $\| \cdot \|$]
When dealing with vector valued functions with values in a
Hilbert space $H$, we use
$\nm \cdot \nm$ for the norm in $H$ induced by the inner product
$\langle \cdot , \cdot \rangle$. We will use the symbol
$\|\cdot \|$ (usually with a subscript) for the norm in the
function space; thus for a vector valued function $f$, the symbol
$\|f\|_{\infty}$ denotes its $L^{\infty}$ norm, which is the
essential supremum of $\nm f(z) \nm$ over $z$ in the domain of
definition of $f$. On the other hand, the symbol $\nm f\nm$ stands
for the scalar valued function whose value at a point $z$ is the
norm of the vector $f(z)$;

\item[$\cdot^{\top}$, $\overline{\;\cdot\;}$, $\cdot^{*}$]
If $M$ is a matrix (possibly infinite), then $M^{\top}$
denotes the transpose of $M$. The complex conjugate of $M$
is denoted by $\overline{M}$, and $M^{*}:=(\overline{M})^{\top}$;

\item[$H^{\infty}$]  space of bounded holomorphic functions
on $\mD$ with the supremum norm;

\item[$H^{p}$]  the Hardy space, i.e.~the space  of analytic functions $f$ on $\mD$
such that
$\|f\|_{p}:=\displaystyle
\sup_{0\leq r<1} \int_{\mT} \nm f(r \zeta)\nm dm(\zeta) <\infty$; we will also use the vector-valued Hardy spaces $H^p(E)$ of functions with values in a Hilbert (or Banach) space $E$;

\item[$A$]  space of bounded holomorphic functions on $\mD$
with a continuous extension to $\mT$ with the supremum norm.

\end{entry}

\newcommand{\D}{\mathbb{D}}

\begin{section}{Introduction}

The paper is devoted to the estimates in the corona problem in some
smooth subalgebras of the algebra $H^\infty$ of bounded analytic
functions in the unit disc $\D$.

There main motivation for studying this problem comes from the idea of
``visibility'' or ``$\delta$-visibility'' of the spectrum, introduced
by N.~Nikolski \cite{Nik99}.

Let us recall the main definitions.  Let $\cA$ be a commutative unital
Banach algebra continuously embedded into the space $C(X)$ of all
continuous functions on a Hausdorff topological space $X$, $\cA\subset
C(X)$. The point evaluations $\delta_{x}$ ($x\in X$) given by
\[
\delta_{x}(f)=f(x), \quad f\in \cA,
\]
are multiplicative linear functionals on $\cA$. Hence if $\cA$
distinguishes points of $X$, then we can identify $X$ with a subset of
the maximal ideal space of $\cA$ (the spectrum $\mathfrak{M}(\cA)$ of
$\cA$), that is, $X\subset \mathfrak{M}(\cA)$.

\begin{df}
\label{df_visible_spectrum}
Let $0<\delta\leq 1$. The spectrum of $\cA$ is said to be
$(\delta,m)$-{\em visible} (from $X$) if there exists a constant
$C(m)$ such that for any vector $f=(f_{1},\dots,f_{m})\in \cA^{m}$
satisfying
\begin{equation}
\label{eqn_corona_condition}
\inf_{x\in X}\sum_{k=1}^{m} |f_{k}(x)|^{2}
\geq \delta^{2}>0
\end{equation}
and the normalizing condition
\[
\|f\|^{2}
:=
\sum_{k=1}^{m} \|f_{k}\|_\cA^{2}
\leq 1,
\]
the Bezout equation
\begin{equation}
\label{eqn_Bezout}
g\cdot f:=\sum_{k=1}^m g_k f_k =e
\end{equation}
has a solution $g=(g_1, \ldots, g_m) \in \cA^{m}$ with
\[
\|g\|= \left(\sum_{k=1}^m \|g_k\|_\cA^2 \right)^{1/2}\leq C(m).
\]
The spectrum is called {\em completely} $\delta$-{\em visible} if it
is $(\delta,m)$-visible for all $m\geq 1$ and the constants $C(m)$ can
be chosen in such a way that $\sup_{m\geq 1} C(m)<\infty$.
\end{df}

This is a norm refinement of the usual corona problem for Banach
algebras, and the motivations for the consideration of this problem
can be found in Nikolski \cite{Nik99}.

The classical corona theorem for the algebra $H^{\infty}$, see
\cite{Car62}, says that if the functions $f_k\in
H^\infty=H^\infty(\D)$ satisfy
\begin{equation}
\label{eqn_first_cor_cond}
1\ge  \sum_{k=1}^{m} |f_{k}(z)|^{2}\geq \delta^{2}>0,
\qquad \forall z\in \D,
\end{equation}
then the Bezout equation
\begin{equation}
\label{bezout_eqn}
\sum_{k=1}^{m} g_{k}f_{k}=1
\end{equation}
has a solution $g_1,g_2, \ldots, g_m$, and moreover the solution
satisfies the estimates
\[
\sum_{k=1}^m |g_k(z)|^2 \le C(\delta, m)^2, \qquad \forall z\in \D.
\]
Later refinements obtained independently by M.~Rosenblum
\cite{Ros80} and V.~Tolokonnikov \cite{Tol80}, got the estimate
independent on $m$ and allowed the case $m=\infty$, see Appendix 3 of
\cite{Nik} for modern treatment.

Note that having estimates that are independent of $m$ in the corona
theorem in fact gives us something slightly more than the complete
$\delta$-visibility of the spectrum of $H^\infty$, since the
normalizing condition in \eqref{eqn_first_cor_cond} is weaker than the
corresponding normalizing condition in Definition
\ref{df_visible_spectrum}.

On the other hand there are many algebras with invisible spectrum.  For
example, for the Wiener algebra $W$ of analytic functions,
\[
f=\sum_{k=0}^\infty \widehat f(k)z^k, \quad
\|f\|_W:=\sum |\widehat f(k)|<\infty,
\]
the Corona Theorem holds trivially, that is, the unit disc $\D$ is
dense in the maximal ideal space $\mathfrak M(W)$, but it is in
general impossible to control the norms of solution of the Bezout
equation: the algebra $W$ is not even $(\delta, 1)$-visible for small
$\delta$.

It is general understanding among experts that the estimates hold for
local norms, and may (generally) fail for non-local norms, for example
for norms given in terms of Fourier coefficients.

In this article, we study the following subalgebras of $H^{\infty}$.
Let us recall that $A$ denotes the \emph{disc algebra} of all bounded
analytic functions continuous up to the boundary, $A=H^\infty\cap
C(\T)$.

\begin{df}
\label{df_algebras}
For a positive integer  $n$ define the following algebras:
\begin{enumerate}
\item $\partial^{-n}H^{\infty}$ is the set of all
analytic functions $f$ defined on $\mathbb{D}$
such that $f$, $f'$, \dots, $f^{(n)}$ belong to $H^{\infty}$.

\item $\partial^{-n}A$ is the set of all analytic
functions $f$ defined on $\mathbb{D}$ such that
$f,f',\dots,f^{(n)}$ belong to the disk algebra $A$.

\item More generally, if $S$ be an open subset of $\mT$,
then $\partial^{-n}A_{S}$ is the set of all analytic functions
$f$ defined on $\mathbb{D}$ such that $f,f', \dots, f^{(n)}$
belong to $A_{S}$, where $A_{S}$ denotes the class of functions
defined on the disk that are holomorphic and bounded in $\mD$
and extend continuously to $S$.
\end{enumerate}
The above spaces are Banach algebras with the norm given by
\[
\|f\|
=
\sum_{j=0}^{n}\frac{1}{j!}\|f^{(j)}\|_{\infty}.
\]
\end{df}
The factor $1/j!$ is chosen so the norm satisfies the
estimate $\|f g\|\le \|f\|\cdot \|g\|$.%
\footnote{In the definition of the Banach algebra it is usually
  required that the norm satisfies the estimate $\|fg\|\le \|f\|\cdot
  \|g\|$. However, in a unital Banach algebra, if one is given the
  norm that only satisfies a weaker inequality $\|fg\|\le C\|f\|\cdot
  \|g\|$ (so the multiplication is continuous), there is a standard
  way to replace the norm by an equivalent one, satisfying the
  inequality with $C=1$. Namely, the new norm of an element $f$ is
  defined as the operator norm of multiplication by $f$. It is an easy
  exercise to show that the new norm is equivalent to the original
  one; one needs the fact that the algebra is unital to get one of the
  estimates.}

For a Hilbert space $H$, one can consider the $H$-valued spaces
$\cA(H)$, where $\cA$ is one of the spaces $\p^{-n}H^\infty$,
$\p^{-n}A$, $\p^{-n}A_S$ defined above. Namely, for an analytic
$H$-valued function $f$ we define its norm as
\begin{equation}
\label{equation_norm_on_A_S_n_page_4}
\|f\|
=
\sum_{j=0}^{n}\frac{1}{j!} \|f^{(j)}\|_{\infty},
\end{equation}
where the norm is understood as the $L^{\infty}$ norm of the
vector-valued function with values in $H$. For example, if $H=\ell^2$
(or $H=\C^m$), then for $f=\{f_k\}_{k=1}^\infty =(f_1, f_2, \ldots,
f_k, \ldots )$,
\[
\|f^{(j)}\|_{\infty}
=
\operatorname{essup}_{z\in \mT} \nm f^{(j)}(z)\nm
=
\operatorname{essup}_{z\in \mT}
\left(\sum_{k}|f_{k}^{(j)}(z)|^{2}\right)^{\frac{1}{2}}.
\]

We prove in the paper that the corona theorem with estimates holds for
all these algebras, and that the estimate does not depend on the
number of functions $f_k$. This fact implies complete
$\delta$-visibility of the spectrum for all $\delta>0$.

One of the motivations for studying these algebras comes from control
theory. Namely, for a system (plant) $G$ with coprime factorization
$G=f_1/f_2$, the construction of a stabilizing feedback is equivalent
to solving the Bezout equation
\[
g_1 f_1 + g_2f_2\equiv 1,
\]
with the stabilizing controller given by $-g_1/g_2$. And assuming
that the original plant $G$ (more precisely, its coprime
factorization) has some smoothness, we want to be able to construct
the stabilizing controller with the same smoothness and to be sure
that the smoothness of this stabilizer is controlled by the smoothness
of $G$.

Before proving the corona theorem with bounds for the subalgebras
of $H^\infty$ introduced above in Definition \ref{df_algebras} ,
we remark that the corona theorem itself (without the estimates)
is trivial for them.  Indeed it is easy to show (see Proposition
\ref{p1.3} below) that the maximal ideal space of our algebras
 (for $n\in \mN$) is the closed unit disk. Then the well known
 equivalence of the density of $X$ in the maximal ideal space
 and the solvability of the Bezout equation \eqref{eqn_Bezout} under
 the assumption \eqref{eqn_corona_condition}
 (with $X=\mD$ in our case) gives the corona theorem for our
algebras.

\begin{proposition}
\label{p1.3}
  Let $\cA$ be one of the algebras $\p^{-n}H^\infty$, $\p^{-n}A$,
  $\p^{-n} A_S$ defined above ($n\ge 1$). The maximal ideal space of
  $\cA$ is the closed unit disk.
\end{proposition}

This proposition is definitely not new. It follows, for example
from \cite[Theorem 6.1]{Tolokonnikov_GDA-1991}. This theorem says,
in particular, that for any algebra of functions $\cA$ satisfying
the property %%
\begin{equation}
\label{GD}
\tag{GD} \text{If } f\in \cA\ \text{and }
\lambda >\|f\|_\infty , \ \lambda\in\C, \  \text{then } (f-\lambda)^{-1}\in \cA,
\end{equation}
 its maximal ideal space coincides with the maximal ideal space of
the $L^\infty$-closure of $\cA$.

The algebras we consider clearly satisfy the condition \eqref{GD},
and $L^\infty$-closure of each algebra is the disc algebra $A$,
whose maximal ideal space coincides with the closed unit disc
$\overline \D$.

For the convenience of the reader we present a (very simple) proof of the above Proposition \ref{p1.3}.
\begin{proof}[Proof of Proposition \ref{p1.3}]
Note that $\partial^{-n}H^{\infty}\subset A$, and so point
evaluation at a fixed $\lambda\in \overline{\mD}$ gives a
multiplicative linear functional on $\partial^{-n}A_{S}$.
We will show that every multiplicative linear functional arises
in this manner.

Let $L$ be a multiplicative linear functional and let $\lambda:=L(z)$
(the value of $L$ on the function $f(z)\equiv z$). Then clearly
$L(f)=f(\lambda)$ for polynomials $f$. We show that for any polynomial
$f$
\begin{equation}
\label{equation_ineq_in_soft_proof}
|L(f)| \leq \|f\|_{\infty}.
\end{equation}
This estimate immediately implies that $|\lambda|\le 1$ (apply
\eqref{equation_ineq_in_soft_proof} to the function $f(z)\equiv z$).
Since $\cA\subset A$, any function $f$ in $\cA$ can be approximated by
polynomials in the $L^\infty$-norm. But
\eqref{equation_ineq_in_soft_proof} implies that $L$ is continuous in
$L^\infty$ norm, so formula \eqref{equation_ineq_in_soft_proof} holds
for all $f\in\cA$. Note that in this reasoning we do not need the
density of polynomials in the norm of $\cA$ (which happens only if
$\cA=\p^{-n}A$).

To prove \eqref{equation_ineq_in_soft_proof} let us notice that if
$f\in \cA$ and $\inf_{z\in \mD}|f(z)|>0$, then $f$ is invertible in
$\cA$. Indeed, since $\cA\subset A$, the condition $\inf_{z\in
  \mD}|f(z)|>0$ implies that $f$ is invertible in $A$.

Differentiating $1/f$ $n$ times we get that all its derivatives up to
the order $n$ are in the algebra $H^\infty$ or $A$ or $A_S$, depending
on the algebra $\cA$ we are considering.

Therefore, if $0\not\in \operatorname{clos}\operatorname{range}(f) =
\operatorname{range}(f)$, then $f$ is invertible in $\cA$, and so $f$
does not belong to any proper ideal of $\cA$. Thus $L(f)\neq 0$ for
any maximal ideal (multiplicative linear functional) $L$. Replacing
$f$ by $f-a$, $a\in \C$, we get that if $a\not\in
\operatorname{range}(f)$, then for any multiplicative linear
functional $L$, $L(f)\neq a$, that is, $L(f)\subset
\operatorname{range}(f)$.  Thus $|L(f)|\leq \|f\|_{\infty}$, and
\eqref{equation_ineq_in_soft_proof} is proved.
\end{proof}

\noindent {\bf Plan of the paper.}
In section \ref{section_corona_H_infty} we prove the corona theorem
with estimates on the norm of the solution for the algebra
$\p^{-n}H^\infty$, see Theorem \ref{theorem_regular_corona_1} below.

This result is stronger than the complete $\delta$-visibility of the
spectrum of $\partial^{-n}H^{\infty}$.

We will use this result to show that the corona theorem with the same
estimates holds for the algebras $\p^{-n} A$ and $\p^{-n}A_S$ as well.
That of course would imply that the spectrum of these algebras is
completely $\delta$-visible for all $\delta>0$.

The estimates for the algebra $\p^{-n}A$ will be obtained from the
estimates for $\p^{-n}H^\infty$ by a simple approximation argument.
The same argument will be used to get the estimates for $\p^{-n}A_S$,
with the essential difference that the construction of the
approximating functions is quite involved in this case: the reasoning
``modulo the approximation'' is very similar to the one for the
$\p^{-n}A$.

Note that the results for $n=0$ are quite known. While we cannot give
the exact reference, the fact that the estimates in the corona theorem
for the disc algebra are the same as the estimates for $H^\infty$ is
known to the specialists.  The estimates in the corona theorem for the
algebra $A_S$ were considered by the first author, \cite{Sas05},
although the equality of these estimates to ones for $H^\infty$ was
not mentioned there.

We should also mention that the Corona Theorem for various
algebras of smooth functions was studied by V.~Tolokonnikov,
\cite{Tolokonnikov_GDA-1991}. In particular, the Corona Theorem
(without estimates) for the algebras considered in our paper
follows from his results, see the remark immediately after
Proposition \ref{p1.3} above. For some algebras of smooth
functions he also obtained the Corona Theorem with estimates.

However, the estimates in the  Corona Theorem for the  algebras we
are considering do not follow from his results. Such estimates,
which are the main goal of the present paper, are completely new.
Also new is the fact that the estimates in all of the algebras we
are considering  are the same (for the same $n$), i.e.~that they
do not depend on continuity properties of the last derivative.

\end{section}

\begin{section}{Estimates in the corona theorem for
$\partial^{-n}H^{\infty}$}
\label{section_corona_H_infty}

\begin{theorem}
\label{theorem_regular_corona_1}
Let $n$ be a nonnegative integer, and let $\cA=\p^{-n}H^\infty$. There
exists a constant $C(\delta,n)$ such that for all $f=(f_1, f_2,
\ldots, f_k, \ldots )\in \cA(\ell^2)$ satisfying
\begin{equation}
\label{equation_1}
0<\delta \leq \nm f(z)\nm_{\ell^2} \;\;\textrm{ for all }z\in \mD,
\end{equation}
and
\begin{equation}
\label{equation_2}
\| f\|_{\cA(\ell^2)} \leq 1,
\end{equation}
there exist $g = (g_1, g_2, \ldots, g_k, \ldots )\in \cA(\ell^2)$ such
that
\begin{equation}
\label{equation_a}
\sum_k g_k(z)f_k(z)=1 \;\; \textrm{ for all }z\in \mD,
\end{equation}
and
\begin{equation}
\label{equation_b}
\|g\|_{\cA(\ell^2)}\leq C(\delta,n).
\end{equation}
\end{theorem}

Note that by considering sequences $f= (f_1, f_2, \ldots, f_n,
\ldots)$ with finitely many non-zero entries, one can get the result
about $m$-tuples as an elementary corollary.

\subsection{Preliminaries for the proof}
We want to introduce a different equivalent norm on the space
$\p^{-n}H^\infty$. Namely, for smooth functions on the circle $\T$ let
us consider the differential operator $D$
\[
(Df)(e^{it}):=-i\frac{d}{dt}f(e^{it}).
\]
Define the space $D^{-n}L^{\infty}:= \{f\in L^{\infty}\;|\; D^{k}f\in
L^{\infty}, k=1, 2, \ldots, n\}$. A natural norm on this class is
given by
\begin{equation}
\label{Norm_DL-1}
\sum_{k=0}^n \|f^{(k)}\|_\infty.
\end{equation}
Of course, one can also define this space for the functions with
values in a Hilbert space $H$ with inner product $\langle \cdot, \cdot
\rangle$, and norm $\nm \cdot \nm$.  For our purposes it is more
convenient to consider a different equivalent norm on $D^{-n}L^\infty$
\begin{equation}
\label{equation_norm_DL}
\|f\|:= \nm \widehat{f}(0) \nm + \|D^{n} f\|_{\infty},
\quad f \in D^{-n}L^{\infty},
\end{equation}
where $\widehat{f}(k)$ ($k\in \mZ$) denotes the $k$th Fourier
coefficient of $f$,
\[
\widehat f(k) = (2\pi)^{-1}\int_\pi^\pi f(e^{it}) e^{-ikt}dt.
\]
To show the equivalence of two norms, let us notice that
for $\zeta \in [0,2\pi)$,
\[
f(e^{i\zeta})
=\frac{1}{2\pi}\displaystyle \int_{\zeta-\pi}^{\zeta+\pi}
[f(e^{i\zeta})-f(e^{i\theta})]d\theta
+\widehat{f}(0).
\]
Since
\[
\nm f(e^{i\zeta})-f(e^{i\theta}) \nm
\le \|D f\|_\infty |\theta -\zeta|,
\]
we get by integrating this estimate
\begin{equation}
\label{equation_norm_DL2}
\|f\|_{\infty}\leq
\frac14\|Df\|_{\infty}+\nm \widehat{f}(0)\nm.
\end{equation}
As $\widehat{Df}(0)=0$, $\|Df\|_{\infty}\leq
\frac14\|D^{2}f\|_{\infty}$.  Proceeding in a similar manner we get
\[
\|D^{k}f\|_{\infty}\leq  4^{k-n} \|D^{n}f \|_{\infty} ,
\quad k \in \{1,\dots , n\},
\quad f\in D^{-n}L^{\infty},
\]
so the norms of all derivatives can be estimated by
$\|D^n f\|_\infty$ and $\nm \widehat f(0)\nm$. Therefore the norms
\eqref{Norm_DL-1} and \eqref{equation_norm_DL} are equivalent.

Now we want to find the predual to $D^{-n}L^\infty$. It is easy to see
that if one writes an appropriate duality, then $D^{-n}L^\infty$ is
dual to $L^1$. Namely, it follows from the standard
$L^{1}$-$L^{\infty}$ duality that any bounded linear functional on
$L^1$ can be represented as
\begin{equation}
\label{equation_duality}
L(f)
=
\langle \widehat{f}(0) ,\widehat{{g}}(0) \rangle
+
\int_{\mT} \langle f, D^{n}{g}\rangle dm,
\qquad f\in L^{1}, \quad
\end{equation}
where $g$ is a function in $D^{-n}L^{\infty}$. Moreover, the norm of
$L$ is comparable with the norm $\|g\|\ci{D^{-n}L^\infty}$.  Indeed,
the functional $L$ can be represented as
\[
L(f)= \int_{\mT} \langle f, F\rangle dm,
\quad f\in L^{1},
\]
where $F\in L^\infty$, $\|F\|_\infty =\|L\|$. Let $D^{-1}$ denote
the integration operator, $D^{-1} e^{int}=\frac{1}{n} e^{int}$, $n\ne
0$. Then $D^{-n}(F-\widehat F(0)) + \widehat F(0)=: g \in
D^{-n}L^\infty$ with the norm $\|g\|_{D^{-n}L^\infty}$ comparable to
$\|F\|_\infty$, which immediately implies the representation
\eqref{equation_duality}.

And finally, it is easy to see that $\p^{-n}H^\infty = H^\infty \cap
D^{-n}L^\infty$ and the norm $\| \fdot\|_{D^{-n}L^\infty}$ is
equivalent to the norm in $\p^{-n}H^\infty$. Indeed, since $D(e^{ikt})
= k e^{ikt}$ we conclude that $Df (z) = z f'(z)$ for analytic
polynomials $f= \sum_{k=0}^N a_k z^k$. Iterating the formula $Df (z) =
z f'(z)$ and using the fact that multiplication by $z$ does not change
the norm in $L^\infty(\T)$ we get the estimate
\[
\|D^k f\|_\infty \le C\sum_{j=1}^k \|f^{(j)}\|_\infty,
\qquad k=1, 2, \ldots, n,
\]
which implies that $\|f\|_{D^{-n}L^\infty} \le C
\|f\|_{\p^{-n}H^\infty}$.

To get the opposite inequality, we iterate the identity $f'(z) =
z^{-1} D F (z)$, and since the multiplication by $z^{-1}$ does not
change the $L^\infty(\T)$ norm we get the estimate
\[
\|f^{(k)}\|_\infty \le C\sum_{j=1}^k \|D^j f\|_\infty,
\qquad k=1, 2, \ldots, n.
\]
Using standard approximation reasoning we get that the norms are
equivalent for functions $f\in \operatorname{Hol}(\overline \D)$,
where $\operatorname{Hol}(\overline \D)$ is the set of all functions
analytic in a neighborhood of the closed disc $\overline \D$.  It is
also easy to see that $\p^{-n} H^\infty \cap
\operatorname{Hol}(\overline \D) = \operatorname{Hol}(\overline \D) =
D^{-n} L^\infty \cap H^\infty \cap
\operatorname{Hol}(\overline \D)$.

Finally, for both $X=\p^{-n} H^\infty$ and $X=D^{-n}L^\infty\cap
H^\infty$ we have that $f\in X$ iff $\sup\{\|f_r\|_X:0\le
r<1\} <\infty$, where $f_r(z) := f(rz)$, and,  moreover $\|f\|_X = \lim_{r\to 1-}
\|f_r\|_X$.

Note that the operator $D$ is symmetric, namely, for smooth $f,g$,
integration by parts or use of the Fourier series representations
yields
\begin{equation}
\label{integration_by_parts}
\int_{\mT} \langle Df, g \rangle dm
=
\int_{\mT} \langle f, Dg \rangle dm.
\end{equation}
Therefore, for smooth functions $f$ the duality
\eqref{equation_duality} can be rewritten as
\begin{equation}
\label{equation_duality1}
L(f)
=
\langle \widehat{f}(0) ,\widehat{{g}}(0) \rangle
+
\int_{\mT} \langle D^{n} f, {g}\rangle dm,
\qquad f\in L^{1}, \quad
\end{equation}

\begin{rem} Given a $\Phi\in C^{\infty}(\overline{\mD})$,
there always exists a $\Psi\in C^{\infty}(\overline{\mD})$ such that
$\overline{\partial}\Psi=\Phi$ on some neighbourhood of
$\overline{\mD}$.  Indeed, let $O$ be open and let
$\overline{\mD}\subset O$.  Let $\alpha \in C_{0}^{\infty}(O)$ be
such that $\alpha =1$ on a neighbourhood of $\overline{\mD}$.
Defining $\Psi$ by
\[
\Psi(z)
=
-\frac{1}{\pi} \iint_{\mR^{2}}
\frac{\alpha(\zeta) \Phi(\zeta)}{\zeta -z} dxdy, \quad z\in \mC,
\]
it can be seen that $\Psi\in C^{\infty}(\mC)$ and $\overline{\partial}
\Psi=\Phi$.
\end{rem}

\subsection{Setting up the $\overline \p$-equation}
We will follow the standard way of setting up the
$\overline\p$-equations to solve the corona problem, as presented for
example in \cite{Nik}. We assume that we are given a column vector $f
= (f_1, f_2, \ldots, f_m, \ldots)^\top$ and we want to find a row
vector $g=(g_1, g_2, \ldots, g_m, \ldots)$ satisfying
\[
g\cdot f = \sum_k g_k f_k \equiv 1.
\]
We will use the standard linear algebra conventions, for example
for a matrix $A$, $A^* = \overline A^\top$. In particular, $f^*$ is a
row vector $f^* = (\overline f_1, \overline f_2, \ldots, \overline
f_m, \ldots)$. Also, for two vectors $f, g\in \ell^2$ we will use the
notation $g\cdot f$ for the ``dot product'', $g\cdot f := g^\top f =
\sum_k g_k f_k$.

As usual, it is sufficient to prove the theorem under the additional
assumption that $f$ is holomorphic in a neighborhood of
$\overline{\mD}$. Let $0<r<1$, and set $f_{r}(z)=f(rz)$, $z\in \mD$.
Then $f_{r}\in \textrm{Hol}(\overline{\mD})$, and we have
$\|f_{r}\|\leq 1$, and $\nm f_{r}(z)\nm \geq \delta$ for all $z\in
\mD$. If the statement of the theorem is true for $f$'s in $
\textrm{Hol}(\overline{\mD})$, then there exists a $g_{r}\in
\textrm{Hol}(\mD)$ such that $g_{r}(z)f_{r}(z)=1$ for all $z\in \mD$,
and $\|g_{r}\|\leq C(\delta)$. If we choose $r_{k}\rightarrow 1$ such
that $g_{r}\rightarrow g$ uniformly on compact subsets of $\mD$
(possible by Montel's theorem), then the $g$ satisfies
(\ref{equation_a}) and (\ref{equation_b}) of the theorem.

We suppose therefore that $f\in \textrm{Hol}(\overline{\mD})$ and
(\ref{equation_1}) holds.

Define the row vector $\varphi$,
\[
\varphi =\frac{f^{*}}{\nm f\nm^{2}} .
\]
Then $\varphi \in C^{\infty}(\overline\mD)$, and $\varphi f \equiv 1$
on a neighbourhood of $\overline\mD$. So $\varphi$ solves the Bezout
equation $\varphi f \equiv 1$, but it is not analytic in $\mD$. Note
that
\[
\overline{\partial} \varphi
=
\frac{(f')^{*}}{\nm f \nm^{2}}-\frac{(f')^{*}f}{\nm f\nm^{4}} f^{*}.
\]
If we find a matrix $\Psi$ solving the $\overline\p$-equation
\[
\overline{\partial} \Psi
= \varphi^{\top} \overline{\partial} \varphi
=: \Phi,
\]
then
\[
g:=\varphi+f^{\top} (\Psi^{\top} -\Psi)
\]
will be analytic in $\mD$, since
\begin{eqnarray*}
\overline{\partial} g
&=&
\overline{\partial} \varphi
+
f^{\top} ( \overline{\partial} \Psi^{\top} -\overline{\partial} \Psi)
\quad \mbox{(since }  \overline{\partial} f =0\mbox{)} \\
&=&
\overline{\partial} \varphi
+
f^{\top} ( (\overline{\partial} \varphi)^{\top}\varphi
-\varphi^{\top}\overline{\partial} \varphi)
\quad \mbox{(using }  \overline{\partial} \Psi
=\varphi^{\top} \overline{\partial} \varphi\mbox{)} \\
&=&
\overline{\partial} \varphi
+
((\overline{\partial} \varphi) f)^{\top}\varphi -\overline{\partial}
\varphi
\quad \mbox{(using }  \varphi f\equiv 1\mbox{)} \\
&=&
((\overline{\partial} \varphi) f)^{\top}\varphi
=
(\overline{\partial} (\varphi f))^{\top}\varphi
\quad \mbox{(since }  \overline{\partial} f =0\mbox{)} \\
&=& 0
\end{eqnarray*}
where the last equality follows from the fact that $\varphi f \equiv
1$. Moreover, since the matrix $\Xi=\Psi - \Psi^\top$ is antisymmetric
($\Xi^\top =-\Xi$), we have $f^\top (\Psi-\Psi^\top) f=0$, so
$gf=\varphi f \equiv 1$.

\subsection{Estimates of the solution of the $\overline\p$-equation
from the boundedness of $L$}

Let us see what we need to get the estimate of the norm of the
solution. Since $D^n (\Xi f) = \sum_{k=0}^n {n\choose k} (D^k \Xi)
D^{n-k}f$, the estimates
\[
\operatornamewithlimits{esssup}_{\zeta\in\T} \nm
\Psi^{(k)}(\zeta)\nm \le C< \infty, \qquad k=1, 2, \ldots, n
\]
where $\nm\cdot\nm$ denotes the operator norm of a matrix, imply
the solution $g$ is in the space $D^{-n}L^\infty(\ell^2)$. Since the
solution $g$ we get is analytic, that is exactly what we need.

\newcommand{\tr}{\operatorname{tr}} \newcommand{\HS}{{\mathfrak{S}_2}}
Since the operator norm of a matrix is dominated by the
Hilbert--Schmidt norm $\nm \cdot\nm_\HS$, it is sufficient to estimate
the Hilbert--Schmidt norms of the derivatives, that is, to estimate
the norm of the solution $\Psi$ in the space $D^{-n}L^\infty(\HS)$.
Note that the space $\HS$ of Hilbert--Schmidt operators (matrices) is
a Hilbert space with the inner product $\langle A,B\rangle_\HS :=
\tr{AB^*}=\tr{B^*A} $, so all the previous discussions about norms and
duality for the space $D^{-n}L^\infty$ do apply here.

We estimate the norm of the solution of the $\overline\p$-equation by
duality. Let $\Psi_0$ be any smooth solution of the
$\overline\p$-equation
\begin{equation}
\label{equation_u}
\overline\p \Psi=\Phi
:=
\varphi^{\top} \overline{\partial} \varphi
=
\frac{\overline{f}}{\nm f\nm^{2}}
\left(
\frac{(f')^{*}}{\nm f\nm^{2}}
-\frac{(f')^{*}f}{\nm f\nm^{4}} f^{*}
\right),
\end{equation}
Define the linear functional $L$ on $H^1_0(\HS):= z H^1(\HS)$,
\[
L(h)
= \int_\T \tr\{(D^{n}h)\Psi_0\} dm
= \int_\T \langle D^n h, \Psi_0^*\rangle_\HS dm.
\]
Note that the above expression is well defined on a dense subspace
of smooth functions in $H^1_0(\HS)$, for example on the subspace
$X_0=H^1_0(\HS)\cap \operatorname{Hol}(\overline D, \HS)$.

If we prove that $L$ is a bounded functional on $H^1_0(\HS)$, it can
be extended by Hahn--Banach Theorem to a bounded functional on the
whole space $L^1(\HS)$. That means, according to our discussions of
duality, see \eqref{equation_duality}, \eqref{equation_duality1}, that
there exists a function $\Psi\in D^{-n}L^\infty$,
$\|\Psi\|_{D^{-n}L^\infty}\asymp \|L\|$, such that
\[
L(h)
= \int_\T \tr\{(D^{n}h)\Psi_0\} dm
= \int_\T \tr\{(D^{n}h)\Psi\} dm
\qquad \forall h\in X_0.
\]
Note that $\widehat h(0)=0$ for $h\in X_0$, so the term
corresponding $\langle \widehat f(0), \widehat g(0)\rangle$ from
\eqref{equation_duality}, \eqref{equation_duality1} disappears.

Since $\int_\T \tr\{(D^{n}h)(\Psi-\Psi_0)\} =0$ on a dense set $X$ in
$H^1_0$, the function $\Psi-\Psi_0$ is analytic in $\D$, so $\Psi $
solves the $\overline\partial$-equation $\overline\p \Psi =\Phi$.

\subsection{Estimates of the functional $L$}
\label{s2.4}
To estimate $L(h)$, we use Green's formula,
\begin{equation}
\label{GrFo}
\tag{G}
\int_\T u\,dm - u(0)
=
\frac{2}{\pi} \iint_\D (\p\overline\p u(z) )
\ln\frac{1}{|z|} dx dy
\end{equation}
which holds for $C^2$-smooth functions $u$ in the closed disc
$\overline\D$ (recall that $\p\overline\p= \frac14 \Delta$). Applying
this formula to $u= \tr\{(D^{n}h)\Psi\}$, where $D^n h$ in the disc is
defined as the harmonic (analytic) extension from the boundary, we get
\begin{eqnarray*}
L(h)
&=&
\int_{\mT} \tr\{ (D^{n}h ){\Psi}\} \, dm
\\
&=&\frac{2}{\pi}\iint_{\mD}
(\partial \overline{\partial} \tr\{ (D^{n} h) {\Psi} \} )
\log \frac{1}{|z|} dx dy  \quad (\text{because }D^nh(0)= 0)
\\
&=&
\frac{2}{\pi} \iint_{\mD} (\partial \tr\{ ( D^{n} h) \Phi \})
\log \frac{1}{|z|} dxdy
\quad \mbox{(because } \overline{\partial}(D^{n}h)=0 \mbox{ and }
\overline{\partial} \Psi=\Phi\mbox{)}
\\
&=&\frac{2}{\pi} ( I_{1}+I_{2} ),
\end{eqnarray*}
where
\[
I_{1}
:=
\displaystyle \iint_{\mD}
\tr\{ (D^{n} h) {\partial} \Phi \}
\log \frac{1}{|z|} dx dy
\quad \text{ and }
I_{2}
:=
\iint_{\mD} \tr\{ (\partial D^{n} h) \Phi \}
\log \frac{1}{|z|} dx dy.
\]

To estimate the integrals $I_1$, $I_2$ we would like to move the
derivatives to $\Phi$. To do this, let us extend the operator $D$ to
the whole disc as follows:
\[
Dw(re^{i\theta})=-i\frac{d}{d\theta} w(re^{i\theta}).
\]
Then $Dz^{n}=nz^{n}$ and $D\overline{z}^{n}=-n\overline{z}^{n}$ for
$n\geq 0$, and so for holomorphic $w$, $ Dw(z)=zw'(z)$ and
$D\overline{w}=-\overline{zw'(z)}$.

Note that if we treat $D^n h$ as the ``extended'' operator $D^n$
applied to the function in the disc, we get the same result as before,
when we defined $D^nh$ in the disc as the harmonic (analytic)
extension from the boundary.

\subsubsection{Estimates of $I_1$}
Using the symmetry of $D$, see \eqref{integration_by_parts}, we get
for $I_1$
\begin{align*}
I_1 = \iint_{\mD}
\tr\{ (D^{n} h) {\partial} \Phi \}
\log \frac{1}{|z|} dx dy
& =\iint_{\mD}
\langle D^{n} h, \overline{\partial} \Phi^* \rangle_{\HS}
\log \frac{1}{|z|} \, dx dy \\
& =
\iint_{\mD} \langle  h, D^{n} \overline{\partial} \Phi^* \rangle_{\HS}
\log \frac{1}{|z|} \, dx dy
\end{align*}
where the last equality can be seen as follows: we write the integral
in polar coordinates, then, in the integral with respect to $d\theta$
we apply the formula (\ref{integration_by_parts}) and finally we go
back to $dx dy$. Note that we used the inner product notation, because
the symmetry of the operator $D$ is more transparent and is easier to
write this way.

Applying $n$ times the operator $D$ we get that $D^n\overline\p
\Phi^*$ can be represented as a sum of terms of form
\begin{equation}
\label{decomp_D^n}
\frac{ \textrm{product of analytic and antianalytic factors}}{\nm f\nm^{2r}}
\end{equation}
where (up to the transpose) the antianalytic factors can be only of
the form $(f^{(j)})^{*}$, the analytic ones can be only of the form
$f^{(l)}$, $j, l=0, 1, \ldots, n+1$.  Moreover, if one looks at the
derivatives of the maximal possible order $k=n+1$, each term of form
\eqref{decomp_D^n} can have at most one factor $f^{(k)}$ and at most
one factor $(f^{(k)})^*$ (it can have both $f^{(k)}$ and
$(f^{(k)})^*$). Indeed, the direct computations show that the function
$\overline\p \Phi^*$ clearly is represented as such a sum, with the
maximal order of each derivative being $1$. Each differentiation $D$
preserves the form, and increases the maximal order of the derivative
at most\footnote{It can be shown by more careful analysis, that no
  cancellation happens, and the maximal order of the derivative
  increases \emph{exactly} by $1$, but we do not need this for the
  proof: we only need that it cannot increase by more than $1$.} by
$1$.

The terms in the decomposition \eqref{decomp_D^n} of
$D^n\overline\p\Phi^*$ containing both factors $f^{(k)}$ and
$(f^{(k)})^*$ of maximal possible order $k=n+1$ can be estimated by
$C\nm f^{(n+1)}\nm_{\ell^2}^2$. Note that $f^{(n)}\in
H^\infty(\ell^2)$.

It is well known (see Section \ref{s2.5} below for all necessary information
about Carleson measures) that for a bounded analytic function $F$ with
values in a Hilbert space the measure $\nm F'(z)\nm^2 \log \frac1{|z|}
\,dxdy$ is Carleson, with the Carleson norm estimated by
$C\|F\|_\infty^2$. Thus we can conclude that the measure $\nm
f^{(n+1)}\nm_{\ell^2} \log\frac1{|z|}\, dxdy$ is Carleson. Therefore
\[
\iint_\D \nm h(z)\nm_\HS \cdot\nm f^{(n+1)}\nm^2_{\ell^2}
\log\frac1{|z|}\, dxdy \le C \|h\|_{H^1(\HS)}
\]
so the terms of $I_1$ containing both $f^{(n+1)}$ and
$(f^{(n+1)})^*$ are estimated.

The terms in the decomposition \eqref{decomp_D^n} of $D^n \overline\p
\Phi^*$ containing only the derivatives of order $k<n+1$ are bounded,
so the corresponding terms in $I_1$ are easily estimated, because the
measure $\log\frac1{|z|}\,dxdy$ is trivially Carleson.

Finally, the terms in \eqref{decomp_D^n} containing only one of the
factors $f^{(n+1)}$ or $(f^{(n+1)})^*$ can be estimated by $C \nm
f^{(n+1)}\nm_{\ell^2}$, and since by the Cauchy--Schwartz inequality
\begin{multline*}
\iint_\D \nm h(z)\nm_\HS \nm f^{(n+1)}(z)\nm_{\ell^2}
\log\frac1{|z|}\, dxdy  \\
\le
\left( \iint_\D \nm h(z)\nm_\HS
\nm f^{(n+1)}(z)\nm_{\ell^2}^2 \log\frac1{|z|}\, dxdy\right)^{1/2}
\left( \iint_\D \nm h(z)\nm_\HS  \log\frac1{|z|}\, dxdy\right)^{1/2}
\\
\le C \|h\|_{H^1(\HS)}
\end{multline*}
(as we discussed above, the measures in both integrals in the second
line are Carleson), so the corresponding terms in $I_1$ are also
easily estimated.\footnote{A careful analysis of
  $D^n\overline\p\Phi^*$ can show that the terms containing only one
  derivative of the maximal order are impossible here, but the above
  reasoning is significantly simpler than the careful analysis of
  derivatives.}

\subsubsection{Estimates of $I_2$}
Let us now estimate $I_2$.  By trivial estimates we have
for $|z|<1/2$,
\[
|\tr\{ (\partial D^{n} h) \Phi \}| \le C\|h\|_{H^1(\HS)},
\]
so we need only to estimate the integral $I_2'$, where one
integrates over $1/2\le |z|<1$.

Indeed, the derivatives of $h$ can be estimated by the standard
estimates for power series, if one recalls that $\nm \widehat h(k)
\nm_{\HS} \le \|h\|_{H^1(\HS)}$. We also have $\nm\Phi(z)\nm \le C\nm
f'(z)\nm$, and using the similar reasoning with power series one can
show that $\nm f'(z)\nm \le C$ for $|z|<1/2$.

Note that for analytic $f$ we have $\partial f=z^{-1}Df$, and so we
can replace $\partial D^{n} h$ by $z^{-1}D^{n+1}h$ in $I_{2}'$. Thus
\begin{equation*}
\label{equation_I2}
I_{2}'
=
\iint\limits_{1/2\le |z|< 1} \tr\{(\partial D^{n}h) \Phi \}
\log \frac{1}{|z|} \, dx dy
=
\iint\limits_{1/2<|z|<1}
\langle z^{-1} D^{n+1}h, \Phi^* \rangle_\HS
\log \frac{1}{|z|}\, dx dy.
\end{equation*}
Using the symmetry of $D$ we get as in the case of $I_1$
\begin{eqnarray*}
I_2'
&=&
\iint\limits_{1/2<|z|<1} \!\!
\langle  D h, D^n ((\overline z)^{-1}\Phi^* )\rangle_\HS
\log \frac{1}{|z|}\, dx dy
\\
&=&
\iint\limits_{1/2<|z|<1} \!\!
\langle  z^{-1} h'(z), D^n ((\overline z)^{-1}\Phi^* )\rangle_\HS
\log \frac{1}{|z|}\, dx dy.
\end{eqnarray*}
Applying the operator $D$ repeatedly to $(\overline z)^{-1}\Phi^*$, we
get the representation of $D^n((\overline z)^{-1}\Phi^* )$ as the sum
of terms of form \eqref{decomp_D^n}, with slight differences. Namely,
the analytic factors, as in the case of $I_1$ can be of the form
$f^{(l)}$, $l=1, 2, \ldots, n+1$, and the antianalytic factors (and
that is the difference with the case of $I_1$) can only be of the form
$(f^{(j)})^*$, $j=1, 2, \ldots, n$ or $(\overline z)^{-\kappa}$,
$\kappa\ge 1$. And again, any term containing the derivative
$f^{(n+1)}$ of the highest possible order can contain it only once.

We notice that $(\overline z)^{-1}\Phi^*$ has such representation with
$n=0$, and each differentiation preserves the form of the
decomposition and increases the maximal possible order of derivatives
$f^{(l)}$ and $(f^{(j)})^*$ by at most $1$.

To estimate $I_2'$, let $h_1$ be a scalar-valued outer function in
$H^2$ such that $|h_1(\zeta)|^2 = \nm h(\zeta)\nm$ a.e.~on $\T$. Then
$h\in H^1(\HS)$ can be represented as $h=h_1 h_2$, where $h_1\in H^2$
(scalar), $h_2\in H^2(\HS)$, and $\|h_1\|_{H^2}^2 =
\|h_2\|_{H^2(\HS)}^2 = \|h\|_{H^1(\HS)}$.

Since $h' = h_1 h_2' + h_1' h_2$, we can estimates the terms of $I_2'$
containing the derivative $f^{(n+1)} $ of the highest possible order
by
\[
\iint_\D \nm h(z)\nm_\HS \nm f^{(n+1)}\nm_{\ell^2} \log\frac1{|z|}
\, dxdy
\le
\iint_\D (|h_1|\cdot \nm h_2'\nm + |h_1'|\cdot \nm h_2\nm)
\nm f^{(n+1)}\nm_{\ell^2} \log\frac1{|z|} \, dxdy.
\]
Since, as we discussed before, when treating $I_1$, the measure
$\nm f^{(n+1)}(z)\nm_{\ell^2} \log\frac1{|z|} \, dxdy$ is Carleson,
with its Carleson norm bounded by $C\|f\|^2_{H^\infty(\ell^2)}$, we get
\begin{align*}
\iint_\D |h_1|\cdot & \nm h_2'\nm \cdot \nm f^{(n+1)}\nm_{\ell^2}
\log\frac1{|z|} \, dxdy
\\
& \le
\left( \iint_\D |h_1|^2 \nm f^{(n+1)}\nm_{\ell^2}^2
\log\frac1{|z|} \, dxdy \right)^{1/2}
\left( \iint_\D \nm h_2'\nm^2  \log\frac1{|z|} \, dxdy \right)^{1/2}
 \\
& \le C \|h_1\|_{H^2} \|h_2\|_{H^2(\HS)} = C \|h\|_{H^1(\HS)};
\end{align*}
here the first integral in the second line is estimated because the
measure is Carleson, and the second integral is simply the
Littlewood--Paley representation of the norm $\|h_2\|_{H^2(\HS)}$. The
integral $\iint_\D |h_1'|\cdot \nm h_2\nm \cdot \nm
f^{(n+1)}\nm_{\ell^2}\log\frac1{|z|} \, dxdy$ is estimated similarly.

The terms in the decomposition \eqref{decomp_D^n} of $D^n((\overline
z)^{-1}\Phi^* )$ which contain only derivatives of order at most $n$
are bounded. Therefore to estimate the rest of $I_2'$ it is sufficient
to estimate $\iint_\D \nm h'\nm \log\frac1{|z|}\,dxdy$. Decomposing as
above $h=h_1 h_2$ and using the fact that the measure
$\log\frac1{|z|}\, dxdy$ is trivially Carleson, we get the estimate
\begin{align*}
\iint_\D |h_1|\cdot & \nm h_2'\nm \log\frac1{|z|} \, dxdy
\\
& \le
\left( \iint_\D |h_1|^2  \log\frac1{|z|} \, dxdy \right)^{1/2}
\left( \iint_\D \nm h_2'\nm^2  \log\frac1{|z|} \, dxdy \right)^{1/2}
 \\
& \le C \|h_1\|_{H^2} \|h_2\|_{H^2(\HS)} = C \|h\|_{H^1(\HS)};
\end{align*}
The integral $\iint_\D |h_1'|\cdot \nm h_2\nm \log\frac1{|z|} \,
dxdy$, and thus the rest of $I_2'$ is estimated similarly.
\hfill\qed

\subsection{Some remarks about Carleson measures}
\label{s2.5} In this subsection we present for the convenience of
the reader some well known facts about the Carleson measures, that
we have used above in Section \ref{section_corona_H_infty}.

Let us recall that a measure $\mu$ in the unit disc $\D$ is
called the \emph{Carleson measure} if the   embedding
$H^2\subset L^2(\mu)$ holds, i.e. if the inequality
\begin{equation}
\label{2.13}
\int_\D |f(z)|^2\,d\mu(z) \le C \|f\|^2_{H^2}, \qquad \forall f\in \D
\end{equation}
holds for some $C<\infty$. The best possible constant $C$
in this inequality is called the \emph{Carleson norm} of the measure $\mu$.

There is a very simple geometric description of the Carleson
measures, cf \cite{Gar81} or any other monograph about $H^p$
spaces. Namely, a measure $\mu$ is Carleson if and only if
$$
\sup_{\xi\in \T\ , r>0} \frac1r \mu\left\{ z\in\D : |z-\xi|<r \right\} <\infty.
$$
Moreover, the above supremum is equivalent (in the sense of
two sided estimate) to the Carleson norm of the measure $\mu$.

However, in this paper we will only use the following simple and well-known
fact about bounded analytic functions and Carleson measures.

\begin{proposition}
\label{p2.3} If $F$ is a bounded analytic function in the unit
disc with values in a Hilbert space, then measure $\mu$, $d\mu(z)
= \log\frac1{|z|} \nm F'(z)\nm^2 \, dxdy$ is Carleson with its
Carleson norm bounded by $C\|F\|^2_\infty$.
\end{proposition}
Note, that this proposition is not true for functions
with values in an arbitrary Banach space.

Note also, that in the scalar case this and even
stronger proposition  is well known and widely used,
see for example the Garnett's book \cite{Gar81}.

There are several ways to prove this proposition, and
it is easier for us to present the proof here and
save the reader a trip to the library, than to give an exact reference.

Probably the simplest way to prove this proposition is
to refer to the so-called Uchiyama
Lemma, cf \cite[Appendix 3, Lemma 6]{Nik}. This lemma
says that if $u\ge 0$ is $C^2$-smooth bounded subharmonic
function (i.e.~$\Delta u\ge 0$) in $\D$, then the measure
$\Delta u(z) \log\frac1{|z|}\,dxdy$ (where $\Delta$ denotes
the Laplacian) is Carleson with Carleson norm estimated by
$2\pi e \|u\|_\infty^2$. Noticing that for an analytic
function $F$ with values in a Hilbert space
$\Delta \nm u(z) \nm^2 = 4\p\pb\nm u(z) \nm^2 =4 \nm u'(z) \nm^2$
we immediately get the proposition with the constant $C=\pi e/2$.

Another, more elementary way to prove the proposition is to use
the Littlewood--Paley formula. Namely, if we apply the Green's
formula (see \eqref{GrFo} in Section \ref{s2.4}) to the function
$u(z) = \nm f(z)\nm^2$, where $f\in H^2(E)$, $E$ is a Hilbert
space, we get the Littlewood--Paley identity
$$
\frac2\pi \iint_\D \nm f'(z)\nm^2 \log\frac1{|z|}\,dxdy = \| f\|^2_{H^2(E)} -\nm f(0)\nm^2.
$$
Thus, if we define the weight $w$ on $\D$ by $w(z) = \frac2\pi \log\frac1{|z|}$, then
$$
\| f'\|_{L^2(w)} \le \|f\|_{H^2},
$$
where $L^2(w) = L^2(E, w)$ is the weighted Lebesgue space of
functions with values in $E$. Applying this estimate to a function
$f$ of form $f=Fg$, $F\in H^\infty(E)$, $g$ is a scalar-valued
function in $H^2$, we get using the triangle inequality
$$
\|F' g \|_{L^2(w)} \le \|F g'\|_{L^2(w)} + \|F g\|_{H^2} \le \|F\|_\infty \|g'\|_{L^2(w)} + \|F\|_\infty \| g\|_{H^2} \le 2 \|F\|_\infty \|g\|_{H^2}.
$$
But this implies that the measure $\frac2\pi \log1{|z|} \,dxdy$
is Carleson with the Carleson norm at most $4 \|F\|_\infty^2$. \hfill\qed

We should also mention that if a measure $\mu$ is Carleson, the
embedding \eqref{2.13} holds (with the same constant) for the
vector-valued $H^2$-spaces $H^2(X)$ with values in
an \emph{arbitrary Banach space} $X$. To see that it is sufficient
to notice that $\nm f(z)\nm \le |h(z)|$ for all $z\in \D$, where
$h$ is the scalar-valued outer function satisfying $|h(\xi)| = \nm f(\xi)\nm$ a.e.~on $\T$.
\end{section}

\begin{section}{Estimates in the corona theorem for other algebras:
preliminaries and the case of $\partial^{-n}A$}
\label{section_integrated_disk_algebra}

\subsection{Continuity of the best estimate}
\label{continuity_of_C}
For a function algebra $\cA$ (one should think about one of the
algebras from Definition \ref{df_algebras}) let $C(\cA, \delta)$,
$\delta>0$ denote the best possible estimate on the norm of the
solution of the Bezout equation,
\[
C(\cA, \delta) := \sup_{f} \inf\{\|g\|_{\cA(\ell^2)} \,|\, g\cdot
f:= \sum_{k} g_k f_k \equiv 1\},
\]
where the supremum is taken over all $f =(f_1, f_2, \ldots, f_m,
\ldots) \in \cA(\ell^2)$, $\|f\|_{\cA(\ell^2)}\le 1$ and such that
\[
\nm f(z)\nm_{\ell^2} := \left(\sum_k |f_k(z)|^2\right)^{1/2} \ge
\delta.
\]

We will show in the rest of the paper that for the function algebras
from the Definition \ref{df_algebras} the constants $C(\cA, \delta)$
coincide,
\[
C(\p^{-n}H^\infty, \delta)
= C(\p^{-n}A, \delta)
= C(\p^{-n}A_S, \delta).
\]

Note that the inequalities
\[
C(\p^{-n}H^\infty, \delta)\le C(\p^{-n}A, \delta), \qquad
C(\p^{-n}H^\infty, \delta) \le C(\p^{-n}A_S, \delta)
\]
are trivial. Indeed, if $f\in \p^{-n}H^\infty(\ell^2)$ and
satisfies the estimates $\|f\|\le 1$, $\nm f(z)\nm \ge\delta$, the
functions $f_r$, $f_r(z)=f(rz)$ are in $\p^{-n}A$ (and in
$\p^{-n}A_S$) and satisfy the same estimates. Therefore, for any
$\varepsilon>0$ one can find solutions $g^r\in \p^{-n}A(\ell^2)$, $g^r
\cdot f_r \equiv 1$, $\|g^r\|\le C(\p^{-n}A, \delta)+\varepsilon$.
Picking a uniformly convergent on compact sets subsequence $g^{r_k}\to
g$, $r_k\to 1-$ (which is possible by Montel's theorem), we get the
$\p^{-n}H^\infty(\ell^2)$ solution $g$, $g\cdot f\equiv 1$, $\|g\|\le
C(\p^{-n}A, \delta)+\varepsilon$. Since $\e$ is arbitrary, we get the
estimate $C(\p^{-n}H^\infty, \delta)\le C(\p^{-n}A, \delta)$. The
estimate for the algebra $\p^{-n}A_S$ is obtained in absolutely the
same way.

Clearly, if $\cA$ is one of the algebras we are considering in the
paper, the functions $\delta\mapsto C(\cA, \delta) $ are
non-increasing. We can say even more:

\begin{lemma}
\label{lm-continuity_of_C}
Let $\cA$ be one of the algebras $\p^{-n} H^\infty$, $\p^{-n}A$,
$\p^{-n}A_S$ ($n\ge0$). Then the function $\delta\mapsto C(\cA,
\delta)$ is continuous on $(0,1)$.
\end{lemma}

Note that this lemma holds for $n=0$, which corresponds to the
case of algebras $H^\infty$, $A$ and $A_S$.

\begin{proof}[Proof of Lemma \ref{lm-continuity_of_C}]
To prove the continuity it is sufficient to only prove uniform right
semi-continuity, that is, that $C(\cA, \delta)=
\lim_{\alpha\to \delta+} C(\cA, \alpha)$ uniformly in
$\delta\in[\delta_0,1)$, for all $\delta_0>0$.
Because $C(\cA, \delta)$ is a non-increasing function of $\delta$,
it will be sufficient to prove only ``$\le$'' estimate
(but still uniformly in $\delta\ge\delta_0$).

Let $f=(f_1, f_2, \ldots, f_m, \ldots)\in \cA(\ell^2)$,
$\|f\|\le 1$, $\nm f(z)\nm_{\ell^2}\ge\delta$ $\forall z\in \D$.

Consider a new vector $\wt f^\gamma$, which is obtained from $f$ by
adding an extra entry $f_0\equiv \gamma$, $\wt f^\gamma =
(\gamma, f_1, f_2, \ldots, f_m, \ldots)$, where $\gamma>0$ is small. Clearly
\[
\nm \wt f^\gamma (z) \nm \ge \sqrt{\delta^2 + \gamma^2}
\qquad \forall z\in \D.
\]
Also,
\[
\| \wt f^\gamma \|_{\cA(\ell^2)}
=
\sqrt{a^2 + \gamma^2} + \|f\|\ci{\cA(\ell^2)} - a
\le
\sqrt{a^2 + \gamma^2} + 1- a,
\]
where $a=\|f\|_{H^\infty(\ell^2)}$. Note that trivially $a\ge\delta$.

The expression $\sqrt{a^2 + \gamma^2} + 1- a$ is a decreasing function
of $a$, so taking into account that $a\ge\delta$ we can estimate
\[
\| \wt f^\gamma \|_{\cA(\ell^2)} \le \sqrt{\delta^2 + \gamma^2} + 1- \delta.
\]
Therefore the $\cA(\ell^2)$ norm of the vector $(\sqrt{\delta^2 +
  \gamma^2} + 1- \delta)^{-1} \wt f^\gamma$ is at most $1$, and we
have
\[
(\sqrt{\delta^2 + \gamma^2} + 1- \delta)^{-1} \nm \wt
f^\gamma(z)\nm_{\ell^2}
\ge \wt
\delta
= \wt \delta(\gamma)
:=
\frac{\sqrt{\delta^2 + \gamma^2}}{ \sqrt{\delta^2+\gamma^2} +1 -\delta}.
\]
Note that $\wt \delta(\gamma)>\delta$ for $\gamma>0$. That can be
checked by noticing that $\wt\delta(0)=\delta$ and that
$\frac{d\wt\delta(\gamma)}{d\gamma}>0$ if $\gamma>0$. Also, trivially,
$\wt\delta(\gamma)\to \delta$ as $\gamma\to 0+$ uniformly in
$\delta\in [\delta_0, 1)$ for all $\delta_0>0$.

Applying the definition of $C(\cA, \delta)$ to the rescaled function
$(\sqrt{\delta^2 + \gamma^2} + 1- \delta)^{-1} \wt f^\gamma$ and then
scaling everything back, we can find a vector $\wt g^\gamma =
(g_0^\gamma, g_1^\gamma, g_2^\gamma, \ldots, g_m^\gamma, \ldots) \in
\cA(\ell^2)$ such that $\wt g^\gamma\cdot \wt f^\gamma\equiv 1$ and
\[
\|\wt g^\gamma\|_{\cA(\ell^2)}
\le
(\sqrt{\delta^2 + \gamma^2} + 1- \delta)^{-1} C(\cA, \wt
\delta(\gamma)) + \gamma
\le
C(\cA, \wt \delta(\gamma)) + \gamma .
\]
Since $C(\cA, \delta)$ is non-increasing, $C(\cA, \wt
\delta(\gamma)) + \gamma \le C(\cA, \delta_0) +1=:M$ for
$\delta\ge\delta_0$, so we have uniform (in $\gamma$ and $\delta
\ge\delta_0$) bound on the norm of $\wt g^\gamma$.

Define $g^\gamma:= (g_1^\gamma, g_2^\gamma, \ldots, g_m^\gamma,
\ldots) \in \cA(\ell^2)$. Since $1= \wt g^\gamma \cdot \wt f^\gamma =
g_0^\gamma \gamma + g^\gamma\cdot f$ and also $\|\gamma
g_0^\gamma\|_\cA \le \gamma \| \wt g^\gamma\|_{\cA(\ell^2)} \le
\gamma\cdot (C(\cA, \delta_0) + 1)=:M\gamma$, we conclude that $\| 1 -
g^\gamma\cdot f \|_{\cA} \le M\gamma\to 0$ as $\gamma\to 0+$.
Therefore for small $\gamma$ the scalar function $g^\gamma\cdot f$ is
invertible in $\cA$ and moreover $\|( g^\gamma\cdot f)^{-1}\|_{\cA}
\le 1/(1-M\gamma) \to 1$ as $\gamma\to 0+$. Then the function
$(g^\gamma \cdot f)^{-1} g^\gamma$ solves the Bezout equation
$(g^\gamma \cdot f)^{-1} g^\gamma \cdot f \equiv 1$, and
\[
\|
(g^\gamma \cdot f)^{-1} g^\gamma
\|_{\cA(\ell^2)} \le (C(\cA, \wt \delta(\gamma)) + \gamma)/(1-M\gamma).
\]

This inequality implies right semi-continuity of $C(\delta)$. Indeed,
since the right side of the equation
\[
\wt \delta
:=
\frac{\sqrt{\delta^2 + \gamma^2}}{ \sqrt{\delta^2+\gamma^2} +1 -\delta}
\]
is an increasing function of $\gamma$, then for $\delta_0\le\delta
\le\wt\delta\le1$ this equation has a unique solution
$\gamma=\gamma(\delta, \wt\delta)$. Moreover, the function
$\gamma(\delta, \wt\delta)$ is clearly continuous (and thus uniformly
continuous) on $\delta_0\le\delta\le\wt\delta\le1$.

Therefore, given $\delta_0>0$ and $\e>0$ one can find $\kappa>0$ such
that for all $\delta$, $\wt\delta$ satisfying $\delta_0\le\delta
\le\wt\delta \le\delta + \kappa$ the inequality $C(\cA, \delta)\le
C(\cA, \wt\delta)+\e$. The inequality $C(\cA, \wt\delta) \le C(\cA,
\delta)$ is trivial because of monotonicity of $C(\cA, \delta)$.
\end{proof}

\subsection{Estimate in the algebra $\p^{-n}A$}
\label{s3.2}

In this section we are going to prove that $C(\p^{-n}A, \delta)=C(\p^{-n}H^\infty,
\delta)$ for $n\ge0$. We only need to prove that $C(\p^{-n}A,
\delta)\le C(\p^{-n}H^\infty, \delta)$, since, as it was discussed
above, the opposite inequality is trivial. Note that here we do not
need the continuity of $C(\cA, \delta)$ proved above in Section
\ref{continuity_of_C}.

Let $f\in (\partial^{-n}A)(\ell^2)$ satisfy
\[
\nm f(z)\nm_{\ell^2} \ge\delta, \qquad \forall z\in \D,
\]
and $\|f\|\le 1$. By the definition of $C(\p^{-n}H^\infty,
\delta)$, for any $\e>0$ there exists $g\in \p^{-n} H^\infty(\ell^2)$
solving the Bezout equation $g\cdot f \equiv 1$ and such that
$\|g\|\le C(\p^{-n}H^\infty, \delta) + \e$.

If $0<r<1$, then $g_{r}\cdot f_{r}\equiv 1$, where $f_{r}(z):=f(rz)$
and $g_{r}(z)=g(rz)$, $z\in \mD$.  So we can write
\[
g_{r}f=g_{r}\cdot f_{r}+g_{r}\cdot (f-f_{r})=1+\alpha_{r},
\]
where $\alpha_{r}:=g_{r}\cdot (f-f_{r})\in \partial^{-n}A$.  Since
$\|f-f_{r}\|\rightarrow 0$ as $r\nearrow 1$ and $\|g_{r}\|\leq \|g\|$,
we can conclude that $\|\alpha_{r}\|\rightarrow 0$ as $r\nearrow 1$.
Thus for $r$ close to $1$, we have that $1+\alpha_{r}$ is invertible
in $\partial^{-n}A$ and $\| (1+\alpha_r)^{-1}\|\to 1$ as $r\nearrow
1$.

Then $(1+\alpha_{r})^{-1} g_{r}f\equiv 1$, and so $(1+\alpha_{r})^{-1}
g_{r}\in \partial^{-n}A$ is a left inverse of $f$. Moreover, since
$\|g_{r}\|\leq \|g\|\leq C(\p^{-n}H^\infty, \delta)+\e$ and
$\|(1+\alpha_{r})^{-1}\|\rightarrow 1$ as $r\nearrow 1$, it follows
that for $r$ sufficiently close to $1$, $\| (1+\alpha_{r})^{-1}
g_{r}\|\le C(\p^{-n}H^\infty, \delta)+ 2\e$.  Therefore $C(\p^{-n}A,
\delta) \le C(\p^{-n}H^\infty, \delta)+ 2\e$. Since $\e>0$ is
arbitrary, we get the desired estimate.

\subsection{Preliminary estimates in the algebra $\p^{-n}A_S$}
In this section we will show that $C(\p^{-n}A_S, \delta) \le 3
C(\p^{-n}H^\infty, \delta)^2$. To get the sharp estimate
$C(\p^{-n}A_S, \delta) \le C(\p^{-n}H^\infty, \delta)$ one needs to
use more delicate reasoning, presented in Section
\ref{section_integrated_A_S} below.

We should emphasize that the reasoning below works only for $n\ge 1$,
that is, that it does not work for the algebra $A_S$.

Let $f\in \p^{-n}A_S(\ell^2)$, $\|f\|\le 1$ satisfy
\[
\nm f(z)\nm_{\ell^2} \ge\delta, \qquad \forall z\in \D.
\]
Let $\e>0$. By the definition of $C(\p^{-n}H^\infty, \delta)$,
there exists $g\in \p^{-n}H^\infty(\ell^2)$ solving the Bezout
equation $g\cdot f \equiv 1$ and such that $\|g\|\le
C(\p^{-n}H^\infty, \delta) + \e$. Then, as before, $g_r\cdot f_r\equiv
1$ for $0<r<1$, where $f_{r}(z):=f(rz)$, $g_{r}(z)=g(rz)$, $z\in \mD$.
We cannot claim that $f_r \to f$ as $r\nearrow 1$ in the norm of
$\p^{-n}H^\infty$, but, since $\p^{-1}H^\infty\subset A$, one can
easily see that the convergence in the weaker norm of
$\p^{-n+1}H^\infty$ takes place (or, equivalently, in the norm of
$\p^{-n+1}A$, which is the same):
\[
\| f_r - f \|_{\p^{-n+1}H^\infty(\ell^2)} \to 0 \qquad \text{as } r\to 0+.
\]
Therefore,
\[
g_{r}f=g_{r}\cdot f_{r}+g_{r}\cdot (f-f_{r})=1+\alpha_{r},
\]
where $\alpha_{r}:=g_{r}\cdot (f-f_{r})\in \partial^{-n}A_S$, and
$\|\alpha_r\|_{\p^{-n+1}H^\infty(\ell^2)}\to 0$ as $r\nearrow 1$. We
can see that $1+\alpha_r\in \p^{-n+1}A$, so $1+\alpha_r$ is invertible
in this algebra $\p^{-n+1}A$ and $\|(1+\alpha_r)^{-1} -
1\|_{\p^{-n+1}A}\to 0$ as $r\nearrow\infty$.

We can show even more, namely that $1+\alpha_r$ is invertible in
$\p^{-n}A_S$ and estimate its norm in this algebra. Namely, let $\f_r=
(1+\alpha_r)^{-2}$. Then clearly $\f_r\in \p^{-n+1}A$ and
$\|\f_r-1\|_{\p^{-n+1}A}\to 0$ as $r\nearrow 1$. Differentiating we
get
\[
((1+\alpha_r)^{-1})' = - (1+\alpha_r)^{-2}\alpha_r' = -\f_r \alpha_r',
\]
so for the $n$th derivative
\[
((1+\alpha_r)^{-1})^{(n)}
=
\sum_{k=0}^{n-1} {n-1\choose k}\f_r^{(k)} \alpha_r^{(n-k)} .
\]
Note that this derivative is continuous on $S$ (because
$\alpha_r\in \p^{-n}A_S$, $\f_r\in \p^{-n+1}A$), so that
$(1+\alpha_r)^{-1}\in \p^{-n}A_S$.  Since
$\|\alpha_r\|_{\p^{-n+1}A}\to 0$ as $r\nearrow 1$,
\[
\left\| \sum_{k=1}^{n-1} {n-1\choose k}
\f_r^{(k)} \alpha_r^{(n-k)} \right\|_\infty \to 0,
\qquad \text{as } r\nearrow 1,
\]
and so
\[
\limsup_{r\to1-} \|((1+\alpha_r)^{-1})^{(n)}\|_\infty \le
\limsup_{r\to1-} \|\f_{r}\|_\infty \|\alpha_r^{(n)}\|_\infty \le
\limsup_{r\to1-} n! \|\alpha_r\|_{\p^{-n}H^\infty}
\]
But it follows from the definition of $\alpha_r$ that
\[
\|\alpha_r\|_{\p^{-n}H^\infty} \le 2 \| g_r\|_{\p^{-n}H^\infty}
\le 2 \| g\|_{\p^{-n}H^\infty}
\le 2(C(\p^{-n}H^\infty, \delta)+ \e).
\]
Using the fact that $\|(1+\alpha_r)^{-1} -1\|_{\p^{-n+1}A}\to 0$ as
$r\nearrow 1$, we can estimate
\begin{align*}
  \limsup_{r\to1-} \|(1+\alpha_r)^{-1}\|_{\p^{-n}H^\infty} &\le
  \limsup_{r\to 1-} \left(1 +
    \frac1{n!}\left\|((1+\alpha_r)^{-1})^{(n)}\right\|_\infty \right)
  \\
  & \le 1 + 2(C(\p^{-n}H^\infty, \delta) +\e) \le 3(C(\p^{-n}H^\infty,
  \delta) +\e)
\end{align*}
Note that the function $(1+\alpha_r)^{-1} g_r$ solves the Bezout
equation, $(1+\alpha_r)^{-1} g_r\cdot f\equiv 1$ and belongs to
$\p^{-n}A_S$. We can estimate the norm
\begin{align*}
\limsup_{r\to1-} \| (1+\alpha_r)^{-1} g_r\|_{\p^{-n}H^\infty}
& \le
\limsup_{r\to1-} \|(1+\alpha_r)^{-1}\|_{\p^{-n}H^\infty}
\|g\|_{\p^{-n}H^\infty}
\\
& \le
3(C(\p^{-n}H^\infty, \delta) +\e) (C(\p^{-n}H^\infty , \delta) +\e).
\end{align*}
Since $\e>0$ is arbitrary we get $C(\p^{-n}A, \delta)
\le 3C(\p^{-n}H^\infty , \delta)^2$. \hfill\qed

\subsection{Remark on the stable rank of the algebras
  $\p^{-n}H^\infty$, $\p^{-n}A$, $\p^{-n}A_S$} Recall that if $R$ is
any ring, then its {\em Bass stable rank}, denoted by
$\textrm{bsr}(R)$, is by definition the least $m$ such that whenever
$r_{1},\dots,r_{m+1} \in R$ and $\{r_{j}\}$ generate $R$ as a left
ideal, there are $b_{1} \dots, b_{m} \in R$ such that
$r_{1}+b_{1}r_{m+1}, \dots, r_{m}+b_{m}r_{m+1}$ generate $R$ as a left
ideal.

The Bass stable rank of each algebra for the function algebras from
the Definition \ref{df_algebras} is equal to $1$. For $n\in \mN$, this
can be deduced easily from the fact that the Bass stable rank of the
disk algebra $A$ is $1$, as follows. (That $\textrm{bsr}(A)=1$ was
shown in \cite{JonMarWol86}.) Suppose that $f_{1},f_{2}\in
\partial^{-n} A_{S}$ generate $\partial^{-n}A_{S}$.  Then
$f_{1},f_{2}\in A$ and for all $z\in \mD$,
$|f_{1}(z)|+|f_{2}(z)|>\delta>0$.  Using $\textrm{bsr}(A)=1$, it
follows that there exists a $g_{2}\in A$ such that $f_{1}+f_{2}g_{2}$
is invertible in $A$. If $r\in (0,1)$, define $g_{2,r}\in
\partial^{-n}A_{S}$ by $g_{2,r}(z):=g_{2}(rz)$, $z\in \mD$. Choosing
$r$ close enough to $1$, we can ensure that $f_{1}+f_{2}g_{2,r}$ is
invertible in $A$, and hence also in $\partial^{-n}A_{S}$.

\end{section}

\begin{section}{Equality of the best estimate in the corona theorem
for $\partial^{-n}A_{S}$ with that for $\p^{-n} H^{\infty}$}
\label{section_integrated_A_S}

In this section we will show that
$C(\p^{-n}A_{S},\delta)=C(\p^{-n}H^{\infty},\delta)$ for $n\ge0$.  The
method is similar to the one used in the previous section for
$\p^{-n}A$, except that we will need a more elaborate approximation
scheme (given in Subsection \ref{subsection_approximation}) below.

The main idea is that we are going approximate the corona data $f$ by
the function $\wt f$ that extends analytically across $S$ to a bigger
(simply connected) domain $\Omega\supset \D$. The solution $\wt g$ of
the Bezout equation $\wt g \cdot \wt f \equiv 1$ restricted to $\D$
automatically belongs to the class $\p^{-n}A_S$ and ``almost solves''
the equation $g\cdot f \equiv 1$. Then, applying the reasoning similar
to the one in Section \ref{s3.2} we get the estimate on the norm of
the solution.

To carry out this plan we first of all need to construct such an
approximation, which is done below in Section
\ref{subsection_approximation}. We will also need to show that we can
keep under control changes of the estimates when we conformally map
$\Omega$ to the disc $\D$.

\subsection{An approximation result}
\label{subsection_approximation}
In this subsection, we prove a result about uniform approximation of a
function from $\p^{-n}A_{S}(\ell^{2})$ by a function holomorphic
across $S$, in Theorem \ref{theorem_approximation}. This result is a
consequence of the following Lemma \ref{lemma_0_approximation}.

\begin{df}
For an open set $\Omega\subset \C$ let $H^{\infty}(\Omega)$ denote
the set of all bounded analytic functions on $\Omega$. If $n$ is a
nonnegative integer, let $\partial^{-n}H^{\infty}(\Omega)$ be the
set of all analytic functions $f$ on $\Omega$ such that
$f,f^{(1)},f^{(2)}, \dots, f^{(n)}$ belong to $H^{\infty}(\Omega)$,
with the norm given by
\[
\| f\|_{\p^{-n}H^\infty(\Omega)} =\sum_{k=0}^n
\frac{1}{k!}\|f^{(k)}\|\ci{H^\infty(\Omega)}.
\]
\end{df}
Note that the space $\p^{-n}H^\infty(\Omega;\ell^2)$ of
$\ell^2$-valued functions is defined similarly. Sometimes, when it is
clear from the context that we are dealing with vector-valued
functions, we will use $\p^{-n}H^{\infty}(\Omega)$ instead of
$\p^{-n}H^\infty(\Omega;\ell^2)$

\begin{lemma}
\label{lemma_0_approximation}
Let $\Omega$ be an open bounded subset of $\C$ containing $0$ and with
boundary $\partial \Omega$ that has a $C^N$-smooth polar
parameterization $r=\rho(\theta)$.  Suppose that $C$ is a closed
subarc in $\partial \Omega$, and $K$ is an open (in $\p\Omega$) set
containing $C$. Let $R$ be the open sector corresponding to $K$, $R=\{
r\zeta: r\ge0,\;\zeta=\rho( \theta) \in K\}$.

Suppose that $f\in \p^{-n}H^\infty(\Omega)=
\p^{-n}H^\infty(\Omega;\ell^{2})$, where $n\le N$, is such that $f$
and all its derivatives $f^{(k)}$ for $k=1, 2, \ldots, n$ extend
continuously to $K=R\cap \p \Omega$.

Then given any $\e>0$, there exists a domain $\wt\Omega= \Omega\cup
O$, where $O$ is an open neighborhood of $C$ in $\mC$ and a
holomorphic function $F:\wt\Omega \rightarrow \ell^{2}$ with the
following properties:
\begin{itemize}
\item[(S1)] $\|F|_{\Omega}-f\|<\e$.

\item[(S2)] The derivatives $F^{(k)}$, $k=0,1, 2, \ldots, n$
extend continuously to $\wt K:=\p\wt\Omega\cap R$.

\item[(S3)] $\left| \|F\|\ci{\p^{-n}H^\infty(\wt \Omega)} -
\|f\|\ci{\p^{-n}H^\infty( \Omega)} \right| <\e$.

\item[(S4)] The boundary $\p\wt \Omega$ of $\wt\Omega$ has a
$C^N$-smooth polar parameterization $r=\wt\rho(\theta)$, and
moreover $\|\rho -\wt\rho\|_{C^N}<\e$.
\end{itemize}
\end{lemma}

\begin{proof}
Define a (trivial radial) $C^{n}$ extension of $f$
(denoted by the same letter) to $\Omega\cup R$ by
\[
f(rz) = f(z), \qquad z\in\p \Omega, \ r>1.
\]
\begin{figure}[h]
   \center
   \psfrag{z}[c][c]{$0$}
   \psfrag{I}[c][c]{$K$}
   \psfrag{V}[c][c]{$R$}
   \psfrag{C}[c][c]{$C$}
   \psfrag{U}[c][c]{$U$}
   \psfrag{W}[c][c]{$W$}
   \psfrag{O}[c][c]{$\Omega$}
   \psfrag{G}[c][c]{$\partial \Omega$}
   \includegraphics[width=6 cm]{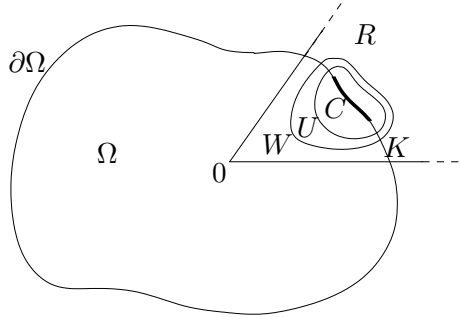}
   \caption{Support of the cut-off function $\varphi$ is contained in $W$.}
   \label{figure_ASn_2}
\end{figure}
Let $\varphi $ be a compactly supported $C^\infty$-function such that
$0\leq \varphi \leq 1$, $\varphi=1$ on a neighbourhood $U$ of $C$ (in
$\mC$), and $\varphi=0$ outside a slightly larger neighbourhood $W$;
see Figure \ref{figure_ASn_2}.

Define a function $h$ (with values in $\ell^2$) by
\begin{equation}
\label{u_12}
h(\zeta)
=
\frac{1}{\pi} \iint
(\overline{\partial} \varphi (z))
\frac{f(z)}{z-\zeta} dx dy
+ \f(\zeta) f(\zeta)
=: u + \f f
\end{equation}
Note that the function $h$ is well-defined for all $z\in \C$, if we
put $\f f = (\overline \p \f) f=0$ outside of $\Omega\cup R$, where
$f$ is not defined.

Moreover $h\in C^n(\C)$. Indeed, the integral $u$ belongs to
$C^{n}(\C)$ since the convolution of the locally integrable function $
z\mapsto \frac{1}{z}$ with the compactly supported $C^n$-function
$(\overline \p \f) f$ is $C^{n}(\mC)$, and trivially $\f f \in
C^n(\C)$.

Using Green's Theorem, one can  see that the formula
\[
u(\zeta) = \frac{1}{2\pi i} \iint \frac{\psi(z)}{\zeta-z} dz\wedge d\overline z=
\frac{1}{\pi} \iint \frac{\psi(z)}{z-\zeta} dxdy
\]
gives, for a continuous compactly supported $\psi$, a solution $u$
of a $\overline\p$-equation $\overline \p u =\psi$; see for instance
\S1 in Chapter VIII of Garnett \cite{Gar81}. Hence, $u$ satisfies the
$\overline{\partial}$-equation
\begin{equation}
\label{d-bar-aux}
\overline{\partial}u=(\overline{\partial}\varphi)f.
\end{equation}

We claim that $h$ is holomorphic in $\Omega$.
Indeed, since $f$ is holomorphic in $\Omega$,
the $\overline\p$-equation \eqref{d-bar-aux} implies
\[
\overline{\partial} h
=
\overline{\partial} (u-\f f)
=
(\overline{\partial}\varphi)f-(\overline{\partial}\varphi)f
=0.
\]
Furthermore, we show that $f-h$ is holomorphic in $U$.  Using again
\eqref{d-bar-aux} and recalling that $\f\equiv 1$ in $U$, we get
$\overline \p u\equiv 0$, $\overline \p \f\equiv 0$ on $ U$, so
$\overline\p h=\overline \p( u -\f f) = \f \overline \p f=\overline\p
f$ in $U$. But that exactly means $f-h$ is analytic in $U$.

\medskip

We observe that if we take the function $F$ to be $f-h$, then it is
holomorphic in $\Omega\cup U$, but it does not necessarily satisfy
condition (S1).  We rectify this situation by adding a shifted version
of $h$ (which is close to $h$).

For $0<r<1$ define $h_{r}(z):=h(rz)$. Since $h\in C^n(\C)$,
\[
h^{(k)}(rz) \to f^{(k)}(z)
\qquad \text{as }r\to 1, \qquad k=0, 1, 2, \ldots, n
\]
uniformly on compact subsets of $\C$. Therefore, we can find $r<1$
sufficiently close to $1$ so that
\begin{equation}
\label{equation_op_cor_3}
\|(h_{r}-h)|_{\Omega}\|\le \e/2 <\e.
\end{equation}
Define $F=f-h+h_{r_{0}}$ on $\Omega\cup R$. The condition (S1) is
satisfied since $\|F-f\|=\|h_{r}-h\|<\e$ on $\Omega$.  Moreover, $F$
is holomorphic in $(\Omega\cup U)\cap \frac1r\Omega =\Omega\cup (U\cap
\frac1r\Omega)=\Omega \cup O_1$ because $f,h,h_{r_{0}}$ are all
holomorphic in $\Omega$, $f-h$ is holomorphic in $U$, and $h_{r}$ is
holomorphic in $\frac1r \Omega$.

Clearly, if $O\Subset O_1$ is an arbitrary open neighborhood of $C$,
then for $\wt\Omega= \Omega \cup O$ the condition (S2) holds (because
$f, h,h_{r}\in C^{n}(O_1)$). The notation $O\Subset O_1$ here means
that $\operatorname{clos}O \subset \operatorname{int}O_1$.

Since $F$ is holomorphic in $ O\Supset C$, for every point $\zeta\in
C$ there exists a neighborhood $V_\zeta\subset O$ of $\zeta$ such that
\[
\sum_{k=0}^n \frac{1}{k!} \nm F(\zeta) -F(z)\nm_{\ell^2}
<\e/3\qquad \forall z\in V_\zeta.
\]
Taking into account \eqref{equation_op_cor_3} we conclude from here
that if we replace $O$ by $\cup_{\zeta\in C} V_\zeta$, then the
condition (S3) will be satisfied.

And it is a trivial exercise to show that we can make $O$ smaller so
that the condition (S4) is satisfied.
\end{proof}

Using the result above, we now prove the following result concerning
uniform holomorphic approximation of functions in
$\partial^{-n}A_{S}(\ell^{2})$. In Lemma \ref{lemma_0_approximation},
we produced an approximate extension of a function across a {\em
  compact} arc, but in the following theorem we construct an
approximate extension across an {\em open} arc.

In order to do this, we decompose the open arc into disjoint open
intervals, and furthermore, we will write each open interval as a
union of closed intervals, and these closed intervals will serve as
the compact arcs of Lemma \ref{lemma_0_approximation}: this lemma will
then be used recursively in order to construct the desired extension.

\begin{theorem}
\label{theorem_approximation} Let $S$ be an open subset of $\mT$,
$n\geq 0$, and $f\in \partial^{-n}A_{S}(\ell^{2})$. Then given any
$\e>0$ and $N\ge n$, there exists a domain $\Omega=\D\cup O$,
where $O$ is an open neighborhood  of $S$ in $\mC$ and a function
$F\in \p^{-n}H^\infty(\Omega; \ell^{2})$ such that
\begin{enumerate}
\item $\|F|_{\mD}-f\|_{\p^{-n}H^\infty}<\e$.

\item $\left| \|F\|\ci{\p^{-n}H^\infty(\Omega)}
- \|f\|\ci{\p^{-n}H^\infty( \D)} \right| <\e$.

\item The boundary $\p\Omega$ has a $C^N$-smooth polar
parametrization $r=\rho(\theta)$, and moreover
$\|\rho-1\|_{C^N}< \e$.
\end{enumerate}
\end{theorem}

\begin{proof} Any open set on $\T$ can be represented as a countable
union of disjoint open intervals (arcs). Each open interval can be
represented as a countable union of \emph{closed} intervals, so we
can represent the open set $S$ as $S=\cup_{n=1}^{\infty} Q_{n} $,
where $Q_{1}$, $Q_{2}$, $Q_{3}$, $\dots$ are closed intervals.

Applying inductively Lemma \ref{lemma_0_approximation} we construct an
increasing sequence of domains $\Omega_{k}$ (in $\mC$) and functions
$\varphi_{k}\in \p^{-n}H^\infty(\Omega_k; \ell^2)$ with the following
properties:
\begin{enumerate}
\item $\Omega_0=\D$, $\f_0=f$.

\item $Q_j\subset \Omega_k$ for $j=1, 2, \ldots, k$.

\item The boundary of $\Omega_k$ has a $C^N$-smooth polar
representation $r=\rho_k(\theta)$, and moreover
$\|\rho_k-\rho_{k-1}\|_{C^N}< \e 2^{-k}$;

\item $\varphi_{k}\in \p^{-n}H^\infty(\Omega_k, \ell^2)$ and
its derivatives $\f^{(j)}$, $j=0, 1, \ldots, n$ extend continuously
to the radial projection $S_k$ of the set $S$ onto $\p\Omega_k$,
$S_k:=\{\rho_k(\theta) e^{i\theta}:\theta \in S\}$.

\item  $\|\varphi_{k}|_{\Omega_{k-1}}-\varphi_{k-1}\|_{\p^{-n}
H^\infty(\Omega_{k-1})} <  {\e}{2^{-k}}$.

\item $\left| \|\f_k\|\ci{\p^{-n}H^\infty(\Omega_k)} -
\|\f_{k-1}\|\ci{\p^{-n}H^\infty( \Omega_{k-1})} \right|
<\e 2^{-k}$
\end{enumerate}
As we mentioned above, we start with $\Omega_0=\D$, $\f_0=f$.  Suppose
$\Omega_{k-1}$, $\f_{k-1}$ are constructed. To get $\Omega_k$, $\f_k$
we apply Lemma \ref{lemma_0_approximation} to the pair $\Omega_{k-1}$,
$\f_{k-1}$ with $2^{-k} \e$ for $\e$. For the arc $C$ we take the
radial projection $C_k$ of $Q_k$ onto $\p\Omega_{k-1}$, $C_k=
\rho_{k-1}(\theta) e^{i \theta}$, and for $K$ the radial projection
$S_{k-1}$ of $S$, $S_{k-1}:= \{\rho_{k-1}(\theta) e^{i\theta}:\theta
\in S\}$.

We need the above assumption (4) to be able to successfully apply
Lemma \ref{lemma_0_approximation}. Condition (4) implies that the
sequence $\f_j$ converges uniformly on each $\Omega_k$, so $F=\lim_j
\f_j$ is an analytic function on $\Omega:= \cup_k \Omega_k$.

Conditions (5) and (6) imply the conclusions (1) and (2) of the
theorem.  Condition (3) on $\f_k$ implies the smoothness of $\p\Omega$
(conclusion (3) of the theorem).
\end{proof}

The above Theorem \ref{theorem_approximation}, for the case $n=0$ and
complex valued functions, can be found in Stray \cite{Str70} and
Gamelin and Garnett \cite{GamGar70}. We will use Theorem
\ref{theorem_approximation} in Subsection \ref{subsection_final}, in
order to prove the estimates in the corona theorem for
$\partial^{-n}A_{S}$.

\subsection{Estimates in the algebra $\p^{-n}H^{\infty}(\Omega)$}
\label{subsection_second_last}

We will prove the corona theorem with bounds for
$\partial^{-n}H^{\infty}(\Omega)$ by using the corresponding result
for $\partial^{-n}H^{\infty}$ obtained earlier, via a conformal map
taking $\mD$ to $\Omega$. We will need the following result by Specht
(see Theorem V and the remark following it, on pages 185--186 of
\cite{Spe51}), which gives bounds on the derivatives of a conformal
map from $\mD$ to $\Omega$, when the boundary of $\Omega$ is smooth
and ``close'' to $\mT$.

\begin{proposition}
\label{proposition_Specht}
Let $C$ be a closed Jordan curve which satisfies the following
assumptions:

\begin{itemize}
\item[(A1)] Every ray from the origin intersects the curve in exactly
one point, and there exists an $\e' \in (0,1)$ such that $C$ lies in
the ring $\{w\in \mC\;|\;1\leq |w|<1+\e'\}$.

\item[(A2)] Let the polar parameterization of $C$ be given by
$\theta (1+\rho(\theta))e^{i\theta}$, $\theta\in[0,2\pi]$,
where $\rho(\theta)$ is nonnegative, and $\rho \in C^{n}$.
Define $\kappa (\theta)=\frac{\rho'(\theta)}{1+\rho(\theta)}$,
$\theta \in [0,2\pi]$. Let $|\kappa'(\theta)|<\e'/\pi$
and $|\omega^{(k)}(\theta)|<\e'/\pi$, $k\in \{2,\dots, n-1\}$,
where $\omega(\theta)=-\arctan (\kappa(\theta))$ (the principal
value of the arctangent is chosen here).

\item[(A3)] For all $\theta_{0}\in [0,2\pi]$,
\[
\frac{1}{2\pi} \int_{-\pi}^{\pi} \left|
\frac{\omega^{(n-1)}(\theta)-\omega^{(n-1)}(\theta_{0})}{\sin
(\frac{\theta-\theta_{0}}{2})} \right|d\theta \leq \e'.
\]
\end{itemize}
Let $\varphi$ be any conformal map $\varphi$ mapping $\mD$ onto the
interior $\Omega$ of $C$ in such a manner that $\varphi(0)=0$ and
$\varphi'(0)>0$. Then $\varphi^{(n)}(z)$ exists for $z \in \textrm{\em
  clos}(\mD)$, and there exist absolute constants $J_{1},\dots,J_{n}$
(that is, numbers which depend only on $n$, but not on $\e'$ or the
curve $C$), such that $|\varphi'(z)-1|\leq J_{1}\e'$ and
$|\varphi^{(k)}(z)|\leq J_{k} \e'$, $k\in \{2,\dots, n\}$.
\end{proposition}

\begin{rem}
\label{remark_Specht} The assumptions (A1)--(A3) of above
proposition are satisfied if $\|\rho\|_{C^{n+1}}<\e $ for
appropriately small $\e$, with $\e'=\e'(\e)$, $\e'(\e)\to 0$ as
$\e\to0$.

The conclusion of the proposition implies the conformal map $\varphi$
belongs to $\p^{-n}A$ is close to the map $z\mapsto z$ in the norm of
$\p^{-n}H^{\infty}$, $\|\f - z\|_{\p^{-n}H^{\infty}}<\gamma(\e')$,
$\gamma(\e')\to 0$ as $\e'\to 0$.
\end{rem}

We now prove the following:

\begin{theorem}
\label{theorem_corona_H_infty_n_Omega} Let $n$ be a nonnegative
integer. Let $\Omega$ be the simply connected open set with boundary a
closed Jordan curve satisfying the assumptions (A1), (A2), (A3) from
Proposition \ref{proposition_Specht}, where $\e'$ is such that
$J_{1}\e'<\frac{1}{2}$. Let $\calA=\p^{-n}H^{\infty}(\Omega)$.

Then for all $f=(f_{1},f_{2},\dots, f_{k}, \dots) \in \calA(\ell^{2})$
satisfying
\[
0<\delta \leq \nm f(z)\nm_{\ell^{2}}\;\;
\textrm{ for all }z\in \Omega,
\quad
\textrm{ and }
\quad
\|f\|_{\calA(\ell^{2})}\leq 1,
\]
there exists a $g=(g_{1},g_{2},\dots, g_{k}, \dots)\in
\calA(\ell^{2})$ such that
\[
\sum_{k} g_{k}(z)f_{k}(z)=1
\textrm{ for all }z\in \Omega, \quad \textrm{ and } \quad
\|g\|_{\calA(\ell^{2})} \leq (1+\alpha(\e'))
C(\p^{-n}H^{\infty},\delta),
\]
where $\alpha(\e')\rightarrow 0$ as $\e'\rightarrow 0$.
\end{theorem}

\begin{proof} Let $\varphi:\mD \rightarrow \Omega$ be a
holomorphic map such that $\varphi(0)=0$ and $\varphi'(0)>0$.
By Proposition \ref{proposition_Specht}, according to
Remark \ref{remark_Specht}, the conformal map $\f$ is
close to the identity map $z$.

Differentiating $f\circ \f$ we get that the $\p^{-n}H^\infty$ norms of
$f$ and $f\circ \f$ are close:
\begin{equation}
\label{equation_for_rescaling_corona_data}
\left|
\|f\|_{\p^{-n}H^{\infty}(\Omega;\ell^{2})} -
\|f\circ \varphi\|_{\p^{-n}H^{\infty}(\D;\ell^{2})}\right|
\leq
\alpha_1\|f\|_{\p^{-n}H^{\infty}(\Omega;\ell^{2})} \le \alpha_1,
\end{equation}
where $\alpha_1=\alpha_1(\e') \rightarrow 0$ as $\e'\rightarrow 0$.
The estimate \eqref{equation_for_rescaling_corona_data} implies that
$\|f\circ\f\|_{\p^{-n}H^{\infty}(\D;\ell^{2})} \le 1+\alpha_1$, so the
``normalized'' vector-function $(1+\alpha_1)^{-1} f\circ\f$ has the
$\p^{-n}H^{\infty}$-norm at most $1$, and satisfies
\[
\frac{1}{1+\alpha_1}\nm f\circ\f(z) \nm_{\ell^2}
\ge
\frac{\delta}{1+\alpha_1}
=:
\wt\delta,
\qquad \forall z\in \D.
\]

Applying to this function the definition of $C(\p^{-n}H^\infty,
\delta)$, we get by solving the Bezout equation for $(1+\alpha_1)^{-1}
f\circ\f$ and then scaling everything back, that there exists
$\widetilde{g}\in \partial^{-n}H^{\infty}(\D;\ell^{2})$ such that
\[
(\wt g\cdot (f\circ\f) )(z) :=
\sum_{k}\widetilde{g}_{k}(z) (f_{k}\circ \varphi)(z)=
1\qquad \forall z\in \mD,
\]
and
\[
\|\widetilde{g}\|_{\p^{-n}H^{\infty}(\D;\ell^{2})}<
(1+\alpha_1)^{-1} C(\p^{-n}H^{\infty},\wt\delta) + \e' \le
C(\p^{-n}H^\infty, \wt\delta) + \e'.
\]

Recalling the continuity of $\delta\mapsto C(\cA,\delta)$, see Lemma
\ref{lm-continuity_of_C}, and noticing that
$\wt\delta=\wt\delta(\e')\to \delta$ as $\e'\to 0$, we can get from
the last estimate that
\[
\|\widetilde{g}\|_{\p^{-n}H^{\infty}(\D;\ell^{2})}<
(1+\alpha_2)
C(\p^{-n}H^{\infty},\delta) .
\]
where $\alpha_2=\alpha_2(\e')\to 0$ as $\e'\to 0$.

Finally defining $g\in \p^{-n}H^\infty(\Omega, \ell^2)$ by
$g:=\widetilde{g} \circ \f^{-1}$ we get the solution of the Bezout
equation $g\cdot f \equiv 1$.  Using
\eqref{equation_for_rescaling_corona_data} again with $g$ replacing
$f$, we can see that the norms of $g$ and $\wt g = g\circ \f$ cannot
differ too much, so we get the desired estimate on the norm of $g$.
\end{proof}

\subsection{Estimates for $\p^{-n}A_{S}$}
\label{subsection_final}

Using Theorem \ref{theorem_approximation} and Theorem
\ref{theorem_corona_H_infty_n_Omega} from the previous two
subsections, we are now ready to prove the estimates in the corona
theorem for $\partial^{-n}A_{S}$.

\begin{theorem}
\label{theorem_corona_A_S_n}
For an open subset $S\subset T$ and $n\ge 0$ we have
$C(\p^{-n}H^{\infty},\delta)=C(\p^{-n}A_{S},\delta)$.
\end{theorem}

\begin{proof}
Let $\calA=\p^{-n}A_{S}$,and let $f=(f_{1},f_{2}, \dots, f_{k},\dots)
\in \calA(\ell^{2})$ satisfy
\[
0<\delta \leq \nm f(z) \nm_{\ell^{2}} \;\; \text{ for all } z\in \mD,
\quad
\text{ and }
\quad
\|f\|_{\calA(\ell^{2})} \leq 1.
\]

Let $\e>0$, be a small number to be specified later. Applying Theorem
\ref{theorem_approximation} (with this $\e$ and $N=n+1$) to the
function $f$ we get a domain $\Omega\supset \D\cup S$ such that its
boundary admits a $C^{n+1}$ polar parameterization
$z=(1+\rho(\theta))e^{i\theta}$, and $\|\rho\|_{C^{n+1}}<\e$. We also
get a function $F\in \p^{-n}H^{\infty}(\Omega;\ell^2)$ such that the
estimates (1) and (2) from the conclusion of Theorem
\ref{theorem_approximation} are satisfied. Estimate (2) implies
that %%
\begin{equation}
\label{4.5} \|F\|_{\p^{-n}H^\infty(\Omega)}\le 1+\e
\end{equation}
and that
\begin{equation}
\label{4.6} \nm F(z)\nm_{\ell^2}\ge \delta -\e \qquad \forall z\in
\Omega
\end{equation}
Let us assume for a moment that $\e$ and $F$ are fixed. Note that if
we make $\Omega$ smaller, the above estimates \eqref{4.5}, \eqref{4.6}
will still hold. Also, if we make $\Omega$ smaller by replacing $\rho$
by $\gamma\rho$, $0<\gamma<1$, the inclusion $\D\cup S\subset \Omega$
will still hold for this smaller $\Omega$.

In light of Remark \ref{remark_Specht}, if we pick sufficiently small
$\gamma$ the boundary of the ``shrunk'' $\Omega$ will satisfy the
assumption assumptions (A1), (A2), (A3) of Proposition
\ref{proposition_Specht}, and, moreover $\e'$ can be made as small as
we want.

Applying Theorem \ref{theorem_corona_H_infty_n_Omega} to the
rescaled function $(1+\e)^{-1}F$ and then scaling everything back
we get that there exists a $\widetilde{g} \in
\partial^{-n}H^{\infty}(\Omega;\ell^{2})$ such that
\[
\wt g \cdot F := \sum_{k}\widetilde{g}_{k}(z)F_{k}(z)=1
\qquad \forall z\in \Omega,
\]
and
\[
\|\widetilde{g}\|_{\p^{-n}H^{\infty}(\Omega)} \leq
(1+\e)^{-1}(1+\alpha(\e'))C(\p^{-n}H^{\infty},\wt\delta),
\]
where $\wt\delta:= (\delta-\e)/(1+\e)$. Since we consider only small
$\e$, we can assume that $\wt\delta\ge\delta/2$. If we make the other
parameter $\e'$ sufficiently small, we get from here the estimate
\[
\|\widetilde{g}\|_{\p^{-n}H^{\infty}(\Omega)} \leq
C(\p^{-n}H^{\infty},\wt\delta).
\]

Define the scalar function $h\in \p^{n} A_S(\D)$ by $h:= \wt g\cdot f$
(both $f$ and $\wt g$ are clearly in $\p^{-n}A_S$). Note that
\[
\|h-1 \|_{\p^{-n}A_S}= \| \wt g \cdot (f-F) \|_{\p^{-n}A_S} \le
\|\wt g  \|_{\p^{-n}A_S} \e \le C\e,
\]
where $C=C(\p^{-n}H^\infty, \delta/2)$. Therefore, for
sufficiently small $\e$, the function $h$ is invertible in
$\p^{-n}A_S$ and
\[
\|h^{-1} \|_{\p^{-n}A_S}\le \frac1{1-C\e}.
\]

The function $g:= h^{-1} \wt g$ clearly belongs to $\p^{-n}A_S$,
solves the Bezout equation $g\cdot f\equiv 1$, and its norm can be
estimated as
\[
\|g \|_{\p^{-n}A_S}\le \|h^{-1} \|_{\p^{-n}A_S} \|\wt g
\|_{\p^{-n}A_S}\le
\frac{C(\p^{-n}H^{\infty},\wt\delta)}{1-C\e}\,,
\]
where recall that $\wt\delta:= (\delta-\e)/(1+\e)$.

Using the continuity of the function $\delta\mapsto C(\p^{-n}H^\infty,
\delta)$, see Lemma \ref{lm-continuity_of_C} above, we get that by
picking sufficiently small $\e$ in the beginning, we can make this
bound as close to $C(\p^{-n}H^{\infty},\delta)$ as we want.
\end{proof}

\end{section}


\begin{thebibliography}{10}

\bibitem{Car62}
L. Carleson.
Interpolations by bounded analytic functions and the corona problem.
{\em Annals of Mathematics. Second Series}, 76:547-559, 1962.

\bibitem{GamGar70}
T.W. Gamelin and J. Garnett.
Uniform approximation to bounded analytic functions.
{\em Revista de la Uni\'on Matem\'atica Argentina}, 25:87-94, 1970.

\bibitem{Gar81}
J.B. Garnett.
{\em Bounded analytic functions}.
Academic Press, 1981.

\bibitem{JonMarWol86}
P.W. Jones, D. Marshall, T.Wolff.
Stable rank of the disc algebra.
{\em Proceeding of the American Mathematical Society},
no. 4, 96:603-604, 1986.

\bibitem{Nik99}
N.K. Nikolski.
In search of the invisible spectrum.
{\em Annales de l'institut Fourier}, no. 6, 49:1925-1998, 1999.

\bibitem{Nik}
%\bysame,
  N.~K. Nikolski\u{\i}.
\emph{{Treatise on the shift operator}}.
Grundlehren der Mathematischen Wissenschaften
[Fundamental Principles of Mathematical Sciences],
vol. 273, Springer-Verlag, Berlin, 1986,
Spectral function theory, With an appendix by
S.V. Hru\v s\v cev [S.V.  Khrushch{\"e}v] and V.V. Peller,
Translated from the Russian by Jaak Peetre.

\bibitem{Ros80}
M. Rosenblum.
A corona theorem for countably many functions.
{\em Integral Equations Operator Theory}, no. 1, 3:125-137, 1980.

\bibitem{Sas05}
A.J. Sasane.
Irrational transfer function classes, coprime factorization and
stabilization. {\em Research Report CDAM-LSE-2005-10},
Center for Discrete and Applicable Mathematics,
London School of Economics, 2005.

\bibitem{Spe51}
E.J. Specht.
Estimates on the mapping function and its derivatives in conformal
mapping of nearly circular regions.
{\em Transactions of the American Mathematical Society}, no. 2,
71:183-196, 1951.

\bibitem{Str70}
A. Stray.
An approximation theorem for subalgebras of $H^{\infty}$.
{\em Pacific Journal of Mathematics},
35:511-515, 1970.

\bibitem{Tol80}
V.A. Tolokonnikov.
{Estimates in Carleson's corona theorem and finitely generated
ideals in the algebra $H^{\infty}$} (Russian).
{\em Akademiya Nauk SSSR. Funktsional. Anal. i
  Prilozhen.}, {\bf 14} (1980), no. 4, 85-86.

\bibitem{Tolokonnikov_GDA-1991}
V.~A. Tolokonnikov, {Generalized {D}ouglas algebras}, \emph{Algebra i Analiz}
  \textbf{3} (1991), no.~2, 231--252,   \emph{Leningrad Math. J.},   {\bf 2}  (1991),  no. 5, 1143--1158
\end{thebibliography}
\end{document}